\documentclass{article}
\usepackage[english]{babel}
\usepackage{amsmath,amssymb,graphicx,stmaryrd}

\catcode`\<=\active \def<{
\fontencoding{T1}\selectfont\symbol{60}\fontencoding{\encodingdefault}}
\catcode`\>=\active \def>{
\fontencoding{T1}\selectfont\symbol{62}\fontencoding{\encodingdefault}}
\newcommand{\assign}{:=}
\newcommand{\backassign}{=:}
\newcommand{\nin}{\not\in}
\newcommand{\textdots}{...}
\newcommand{\tmem}[1]{{\em #1\/}}
\newcommand{\tmop}[1]{\ensuremath{\operatorname{#1}}}
\newcommand{\tmrsup}[1]{\textsuperscript{#1}}
\newcommand{\tmstrong}[1]{\textbf{#1}}
\newcommand{\tmtextit}[1]{{\itshape{#1}}}
\newcommand{\tmtexttt}[1]{{\ttfamily{#1}}}
\newcommand{\tmverbatim}[1]{{\ttfamily{#1}}}
\newenvironment{itemizedot}{\begin{itemize} }{\end{itemize}}
\newtheorem{corollary}{Corollary}
\newtheorem{lemma}{Lemma}
\newtheorem{theorem}{Theorem}

\begin{document}

\title{The Hardy space from an engineer's perspective}

\author{
  Nicola Arcozzi
  \and
  Richard Rochberg
}

\date{September 25, 2020}

\maketitle

\begin{abstract}
  We give an overview of parts of the theory of Hardy spaces from the
  viewpoint of signals and systems theory. There are books on this topic,
  which dates back to Bode, Nyquist, and Wiener, and that eventually led to
  the developement of $H^{\infty}$ optimal control. Our modest goal here is
  giving a beginner's dictionary for mathematicians and engineers who know
  little of either systems or $H^2$ spaces.
\end{abstract}

\tableofcontents

\section{Introduction}

The theory of Hardy spaces is a nice example of the ``unreasonable
effectiveness of mathematics'' in providing a conceptual and computational
framework for the applied sciences. The theory itself lives comfortably in
pure mathematics. It had its inception in Privalov's study of the boundary
behavior of bounded holomorphic functions, some years before Hardy defined the
spaces which go under its name. For many years the Hardy spaces $H^p$ and the
operators acting on them were studied in great depth, and an elegant and
profound theory was developed.

A notable breakthrough was C. Fefferman discovery, in 1971, that the dual of
the Hardy space $H^1$ is the space $\tmop{BMO}$ of functions having bounded
mean oscillations. This result contained the definite solution of the problem
of characterizing the symbols for which the corresponding Hankel operator is
bounded on $H^2$, developing a line of investigation in which Nehari had been
a primary figure. One of the unxpected features of Fefferman result is that
$\tmop{BMO}$ had been earlier defined by Fritz John, and developed by him and
Luis Niremberg, in the distant realm of elasticity theory (``the unreasonable
effectiveness of mathematics'' of the applied sort in providing tools for the
pure ones).

While the pure mathematicians were developing the theory of the Hardy spaces,
engineers found out that they were a very useful tool in signal processing,
then in linear control theory. The basic idea is that signals and systems can
be extended, in frequency space, to holomorphic functions, whose poles and
zeros provide crucial information. This was the beginning of $H^2$ control
theory. The use of frequency methods was pioneered by Bode, Black, and Nyquist
at Bell Labs in the 1930's. Soon after, Wiener entered the picture designing
optimal filtering. Helton, Francis, and many others, developed the
contemporary theory and applications between 1970's and 1990's.

Our goal here is providing an overview of some rather classical parts of Hardy
space theory, highlighting the interpretation in terms of signals and systems.
We hope this helps the pure mathematician, especially the one who is new to
the topic, to develop an intuition for it. Partial as they are, intuitions are
a necessary part of understanding. On the other side, we aim at convincing the
engineer eventually reading these notes that there are interesting things in
Hardy theory to be learned, interpreted, used.

The frontier between these theories is so vast that we do not even try to make
a list of what we are not covering. For the topics we do cover we will not
give specic references to the literature. We do however include at the end a
list of some the many books and surveys in the area, with the hope they will
help the interested reader who wants to learn more. We restrict to signals in
discrete time. The case of continuous time is not much different, but for
technical headaches. We do not even mention the matrix valued case, that is,
what we say concers SISO (single input/single output) systems, not SIMO or
MIMO ones.

The Hardy space theory functions as a model for those studying holomorphic
function spaces, and often the first questions asked when studying a different
function space are ``do things work here as in the Hardy space?'' In the final
section we discuss that question and others for closely related function
spaces, including the Dirichlet space.

\section{Linear systems without holomorphic functions}

We will work all along with complex valued signals in discrete time, i.e.
$\phi : \mathbb{Z} \rightarrow \mathbb{C}$, the space of which is denoted by
$\ell (\mathbb{Z})$. It will be soon clear that the complex field is best
suited for dealing with linear systems, and real valued signals can be
treated, with some care, as a special case. In doing preliminary calculations
we consider signals $\phi$ with finite support, $\phi (n) = 0$ for $| n |$
large, and write $\phi \in \ell_c (\mathbb{Z})$. A \tmverbatim{single
input/single output system} (SISO) is simply a map $T : \ell (\mathbb{Z})
\rightarrow \ell (\mathbb{Z})$, defined on some subset of allowable signals.

Some properties a system is often required to satisfy are the following.
\begin{itemizedot}
  \item \tmverbatim{Linearity}: $T (a \phi + b \psi) = aT (\phi) + bT (\psi)$,
  in which cas we write $T (\phi) = T \phi$;
  
  \item \tmverbatim{Time (or shift) invariance}: let $\tau_1 \phi (n) = \phi
  (n - 1)$ be the forward \tmverbatim{shift} by one unit of time, then $T
  (\tau_1 \phi) = \tau_1 (T (\phi))$;
  
  \item \tmverbatim{Causality}: if $\phi (n) = \psi (n)$ for all $n \leqslant
  m$, then $T (\phi) (m) = T (\phi) (m)$;
  
  \item $p$-\tmverbatim{Stability}: for a linear system, it can be phrased as
  $\interleave T \interleave_{\mathcal{B}(\ell^p)} = \sup_{\phi}  \frac{\|T \phi \|_{\ell^p}}{\|
  \phi \|_{\ell^p}} < \infty$, where
  \[ \| \phi \|_{\ell^p} = \left\{ \begin{array}{l}
       \sup_n | \phi (n) | \text{if } p = \infty\\
       \left( \sum_n | \phi (n) |^p \right)^{1 / p}  \text{if } 1 \leqslant p
       < \infty
     \end{array} \right. \]
  is a measure of the size of the signal, the choices $p = 1, 2, \infty$ being
  the most important in applications.
\end{itemizedot}
The meaning of time invariance is clear: the system works the same way all
times; if the input $\phi$ is delayed by one time unit, $\tau_1 \phi$, then
the output $T (\phi)$ is delayed by one unit of time. Causality means that the
output $T (\phi) (m)$ at time $m$ only depends on inputs up to time $m$, not
on future information. In other words, the time scale for input and output is
the same: if we process a signal in its entirety, as it is done for instance
when denoising an old musical record, causality is not an issue; but if we
denoise a broadcast in real time, then causality is an obvious requirement.

Stability is a requirement of systems (bounds on energy, on size,{\textdots}),
or, often, a law of nature, if the system describes a phenomenon. The
assumption of linearity simplifies the mathematics and is a very good
approximation to many systems of interest. We will not consider the nonlinear
theory here.

It is an easy and instructive exercise using the definitions to show that a
linear, time invariant system is causal if and only if $\phi (n) = 0$ for
negative $n$ implies $T \phi (n) = 0$ for negative $n$. We will denote by
$\ell (\mathbb{N})$ the subspace of those $\phi$ in $\ell (\mathbb{Z})$ for
which $\phi (n) = 0$ for negative $n$ and we set $\ell_c (\mathbb{N}) = \ell_c
(\mathbb{Z}) \cap \ell (\mathbb{N})$. Causality can then be rephrased as
saying that $T : \ell (\mathbb{N}) \rightarrow \ell (\mathbb{N})$.

The characterization of linear, shift invariant systems acting on $\ell_c
(\mathbb{Z})$, is purely algebraic, as it is that of the subclass of causal
ones. We recall that the \tmverbatim{convolution} of $\phi, \psi : \mathbb{Z}
\rightarrow \mathbb{C}$ is $\phi \ast \psi : \mathbb{Z} \rightarrow
\mathbb{C}$,
\[ \phi \ast \psi (m) = \sum_n \phi (m - n) \psi (n) = \psi \ast \phi (m), \]
whenever the sum is defined (e.g. if $\phi$ or $\psi$ belong to $\ell_c
(\mathbb{Z})$).

\begin{theorem}
  Let $T$ be a linear system defined on $\ell_c (\mathbb{Z})$. Then, $T$ is
  shift invariant is and only if there is a function $k : \mathbb{Z}
  \rightarrow \mathbb{C}$ such that
  \[ T \phi = k \ast \phi . \]
  Moreover $k$, the \tmverbatim{unit impulse response}, is uniquely determined
  by $k = T \delta_0$, where $\delta_m (n) = \left\{ \begin{array}{l}
    1 \text{ if } n = m\\
    0 \text{ if } n \neq m
  \end{array} \right.$ .\quad  The system is also causal if and only if
  \[ k (n) = 0 \text{for } n < 0. \]
\end{theorem}

Let $\tau_m \phi (n) = \phi (n - m) = \tau_1^{\circ m} \phi (n)$, $m = \sigma
|m| \in \mathbb{Z}$, where $f^{\circ m} = f^{\sigma} \circ \ldots \circ
f^{\sigma}$, $|m|$ times. In particular, $\tau_m \delta_n = \delta_{n + m}$.
Then, using time invariance of $T$ in the third equality,
\begin{eqnarray*}
  T \phi (n) & = & T \left( \sum_m \phi (m) \delta_m \right) (n) .\\
  & = & \sum_m \phi (m) T (\tau_m \delta_0) (n)\\
  & = & \sum_m \phi (m) \tau_m T (\delta_0) (n)\\
  & = & \sum_m \phi (m) T (\delta_0)  (n - m)\\
  & = & \phi \ast T (\delta_0) (n) .
\end{eqnarray*}
That the system $\phi \mapsto k \ast \phi$ is time invariant is easy to check.
If $T$ is also causal, then
\[ k (m) = T \delta_0 (m) = 0 \text{for all } m < 0 \]
because $\delta_0 (m) = 0$ for negative $m$.

In the causal case, the action of $T$ on $\phi \in \ell (\mathbb{N})$ is a
finite sum:
\[ k \ast \phi (m) = \sum_{n = 0}^m k (m - n) \phi (n) . \]
Although the algebraic analysis is straightforward, the analytic details are
subtle. The problem lies in establishing stability. We consider here the case
$p = 2$, which will take us to the Hardy spaces, but we first mention $p =
\infty$, leading to Wiener's algebra.

For a linear system (operator) $T : X \rightarrow Y$ between two Banach
function spaces $X$ and $Y$ we write
\[ \interleave T \interleave_{\mathcal{B} (X, Y)} = \sup_{v \in X, v \neq 0}
   \frac{\| T v \|_X}{\| v \|_Y}, \]
and we shorten $\mathcal{B} (X, X) =\mathcal{B} (X)$.

\begin{theorem}
  A linear, time invariant system is $\infty$-stable if and only if $k \in
  \ell^1 (\mathbb{Z})$, in which case $\interleave T \interleave_{\mathcal{B}
  (\ell^{\infty} (\mathbb{Z}))} = \|k\|_{\ell^1}$.
\end{theorem}

The elementary estimate
\[ |k \ast \phi (n) | \leqslant \|k\|_{\ell^1} \cdot \| \phi
   \|_{\ell^{\infty}} \]
gives us $\interleave T \interleave_{\infty} \leqslant \|k\|_{\ell^1}$. In the
other direction, set $\phi (n) = \frac{\overline{k (- n)}}{|k (- n) |} \chi (n
: k (- n) \neq 0)$ to have $k \ast \phi (0) = \|k\|_{\ell^1}$ and $\| \phi
\|_{\ell^{\infty}} = 1$.

We leave it to the reader to show that in the causal case $k \in \ell^1
(\mathbb{N})$, we could consider an extremal sequence $\phi_m \in
\ell^{\infty} (\mathbb{N})$ to show that
\[ \sup_{\phi \in \ell^{\infty} (\mathbb{N})}  \frac{\|k \ast \phi
   \|_{\ell^{\infty} (\mathbb{N})}}{\| \phi \|_{\ell^{\infty} (\mathbb{N})}} =
   \|k\|_{\ell^1 (\mathbb{N})}, \]
i.e. that the $\infty$-norm of a causal system can be estimated by considering
signals in positive time.

The space $\ell^1 (\mathbb{Z)}$ with the multiplication given by convolution
is a Banach algebra. Using Fourier series the algebra is isomorphic to the
Banach algebra of continuous functions on the circle which have absolutely
convergent Fourier series, now with multiplication given by the pointwise
product of functions. Both versions are called the Wiener algebra.

\

The case of $2$-stability is richer.

\begin{theorem}
  We have $\interleave T \interleave_{\mathcal{B} (\ell^2 (\mathbb{Z}))}
  \leqslant \|k\|_{\ell^1}$, with equality if $k \geqslant 0$.
  
  In the causal case, we have
  \[ \||T|\|_{\mathcal{B} (\ell^2 (\mathbb{Z}))} = \sup_{\phi \in \ell^2
     (\mathbb{N})}  \frac{\|T \phi \|_{\ell^2 (\mathbb{N})}}{\| \phi
     \|_{\ell^2 (\mathbb{N})}} . \]
  However there are systems, even stable ones, for which $\|k\|_{\ell^1} =
  \infty$.
\end{theorem}

The estimate $\interleave T \interleave_{\mathcal{B} (\ell^2 (\mathbb{Z}))}
\leqslant \|k\|_{\ell^1}$ follows from an easy instance of Hausdorff-Young's
inequality,
\[ \|k \ast \phi \|_{\ell^p} \leqslant \|k\|_{\ell^1} \cdot \| \phi
   \|_{\ell^p}, \]
which holds for $1 \leqslant p \leqslant \infty$. If $T$ is causal, to have
its norm we can just test on $\phi \in \ell^2 (\mathbb{N})$; this will be
easily proved using holomorphic functions. Using holomorphic theory, examples
with $\interleave T \interleave_{\mathcal{B} (\ell^2 (\mathbb{N}))} < \infty$
and $\|k\|_{\ell^1} = \infty$ will naturally come to mind. Using that approach
we will find necessary and sufficient conditions on $k$ for $T$ to be stable.

A reasonable problem is {\tmstrong{designing a causal system}} $T$,
{\tmstrong{that is as close as possible to a given non causal system}} $V$:
$V$ is what we would like to do, while $T$ is what we can do remaining in the
causal class. A quantitative way to state the problem is the following. For
given $V$ with $\interleave V \interleave_{\mathcal{B} (\ell^2 (\mathbb{Z}))}
< \infty$, we want find a causal $T$ for which it is achieved
\[ \min_{T \text{causal}} \sup_{\phi \in \ell^2 (\mathbb{N})}  \frac{\|V \phi
   - T \phi \|_{\ell^2}}{\| \phi \|_{\ell^2}} . \]
We will see later that the problem has a solution within Nehari's theory of
Hankel operators, which will be sketched below.

Another important problem is having {\tmstrong{the complete library of
time-invariant features of signals}}; that is, those features which remain
unchanged if the signal is anticipated or delayed. One such quality is the
frequency spectrum, which we will more rigorously define below.

Each feature might be identified with the set $\mathcal{H} \subseteq \ell^2$
of the functions $\phi$ having that feature. The time invariance of the
feature can be meant in a strong sense (\tmverbatim{bi-invariance}):
\[ \phi \in \mathcal{H} \Leftrightarrow \tau_1 \phi \in \mathcal{H}, \]
or in a weaker sense ([\tmverbatim{forward}] \tmverbatim{invariance}):
\[ \phi \in \mathcal{H} \Rightarrow \tau_1 \phi \in \mathcal{H}, \]
in which a signal might acquire a feature it did not possess before. This is
especially meaningful in the causal case, where the only bi-invariant (linear)
features are trivial: all or none.

As we are dealing with linear theory, we will assume that $\mathcal{H}$ is a
closed, linear subspace of $\ell^2$, and that $\mathcal{H} \neq 0, \ell^2$ is
not trivial. We will say in this case that $\mathcal{H}$ is a
\tmverbatim{bi-invariant}, resp. \tmverbatim{invariant}, subspace of $\ell^2$.

\section{Time and Frequency, and the time-invariant case.}

In this section we review the $L^2$ Fourier theory on $\mathbb{Z}$, which
might be read as Fourier series upside-down. The first motivation comes from
invariant subspaces. Suppose $\phi \neq 0$ is an eigenfunction of the shift,
$\tau_1 \phi = \lambda \phi$ (with, by necessity, $\lambda \neq 0$). Then,
span$\{\phi\}$ is a $1$-dimensional bi-invariant subspace, provided that $\phi
\in \ell^2$.

A little calculation gives
\begin{eqnarray*}
  \phi (n) & = & \lambda^{- 1} \tau_1 \phi (n) = \lambda^{- 1} \phi (n - 1)\\
  & = & \lambda^{- 1} \tau_1 \phi (n - 2)\\
  & \ldots & \\
  & = & \lambda^{- n} \phi (0),
\end{eqnarray*}
a formula which hold for negative $n$'s as well. After normalizing $\phi (0) =
1$, we see that (i) $\phi \nin \ell^2 (\mathbb{Z})$, and (ii) $\phi$ is
bounded if and only if $\lambda = e^{it}$ for some $t \in (0, 2 \pi]
=\mathbb{T}$, in which case $\phi (n) = e_t (n) = e^{- nit}$. It is natural to
assign to the signal $e_t$ the \tmverbatim{period} $2 \pi / t \geqslant 1$: a
time interval which is a fortiori larger than the gap between successive
integers; then a \tmverbatim{frequency} $\omega = t / 2 \pi$.

To each signal $\phi \in \ell^2$ assign its \tmtexttt{Fourier transform}
$\hat{\phi} (e^{it}) = \sum_n \phi (n) e^{int}$, a function in $L^2 =$ $L^2
(\mathbb{T}, d \theta / 2 \pi)$ with $\| \phi \|_{\ell^2} = \| \hat{\phi}
\|_{L^2}$. Then,
\begin{eqnarray*}
  \| \phi \|_{\ell^2}^2 & = & \int_{\mathbb{T}} | \hat{\phi} (e^{it}) |^2
  \frac{dt}{2 \pi},\\
  \phi (n) & = & \frac{1}{2 \pi}  \int_{\mathbb{T}} \hat{\phi} (e^{it}) e^{-
  int} dt,\\
  (\phi \ast \psi)^{\widehat{}} (t) & = & \hat{\phi} (t)  \hat{\psi} (t) .
\end{eqnarray*}
This is all we need from Fourier theory.

\subsection{The characterization of time invariant operators}

From these relations, it is easy to characterize time invariant operators on
$\ell^2 (\mathbb{Z})$.

\begin{theorem}
  The time-invariant system $T \phi = k \ast \phi$ is $2$-stable if and only
  if $\hat{k} = b \in L^{\infty} (\mathbb{T})$. Moreover,
  \[ \interleave T \interleave_{\mathcal{B} (\ell^2 (\mathbb{Z}))} = \sup
     \frac{\|bh\|_{L^2 (\mathbb{T})}}{\|h\|_{L^2 (\mathbb{T})}} . \]
\end{theorem}

Denote by $M_b : h \mapsto bh$ the operator of multiplication times $b$. Then,
$\interleave T \interleave_{\mathcal{B} (\ell^2 (\mathbb{Z}))} = \interleave
M_b \interleave_{\mathcal{B} (L^2 (\mathbb{T}))}$, where the latter refers to
the norm as bounded operator on $L^2 (\mathbb{T})$.

The proof is easy. First, $k = T \delta_0$ is a priori in $\ell^2$, hence $b$
is in $L^2$, and
\begin{eqnarray*}
  \sum_n |k \ast \phi (n) |^2 & = & \frac{1}{2 \pi}  \int_{\mathbb{T}} | (k
  \ast \phi)^{\widehat{}} (t) |^2 dt\\
  & = & \frac{1}{2 \pi}  \int_{\mathbb{T}} |b (t) \hat{\phi} (t) |^2 dt\\
  & \leqslant & \|b\|_{L^{\infty}}^2 \frac{1}{2 \pi}  \int_{\mathbb{T}} |
  \hat{\phi} (t) |^2 dt\\
  & = & \|b\|_{L^{\infty}}^2 \| \phi \|_{L^2}^2,
\end{eqnarray*}
hence $\interleave T \interleave_{\mathcal{B} (\ell^2 (\mathbb{Z}))} \leqslant
\|b\|_{L^{\infty} (\mathbb{T})}$, and choosing $\hat{\phi} (t)$ supported
where $|b (t) |$ is close to its supremum it is easy to show that $\interleave
T \interleave_{\mathcal{B} (\ell^2 (\mathbb{Z}))} \geqslant \|b\|_{L^{\infty}
(\mathbb{T}) - \epsilon}$ for all positive $\epsilon$.

The function $b = \hat{k}$ is the \tmverbatim{transfer function} of the system
$T \phi = k \ast \phi$.

\subsection{The characterization of bi-invariant and invariant spaces for the
shift on $\ell^2 (\mathbb{Z})$}

Similarly simple is the characterization of the bi-invariant subspaces:
{\tmstrong{the invariant features are the sets of frequencies}}. First, on the
frequency side we look for subspaces $\hat{\mathcal{H}}$ of $L^2 (\mathbb{T})$
such that $S \hat{\mathcal{H}} = \hat{\mathcal{H}}$, where $Sh (t) = e^{it} h
(t)$ is the shift on the frequency side. We still call them ``invariant
subspaces for the shift''.

\begin{theorem}
  $M$ is a closed doubly invariant subspace of $L^2 = L^2 (\mathbb{T})$ if and
  only if $M = \eta L^2$ for some $\eta$ which is the characteristic function
  of some $E \subset \mathbb{T}$.
\end{theorem}

That $M$ is doubly invariant is straightforward.

Suppose we have such an $M$. Let $P$ be the orthogonal projection of $L^2$
onto $M$ and let $\eta = P (1)$. Let $\gamma (t) = e^{it}$. By definition of
the projection $1 - \eta \perp M$, hence $1 - \eta \perp \eta \gamma^n$ for
all $n \in \mathbb{Z}$.
\[ 0 = < 1 - \eta, \eta \gamma^n > = \frac{1}{2 \pi}  \int_{\mathbb{T}}
   (\bar{\eta} - | \eta |^2) \gamma^n dt, \]
so all the Fourier coefficients of $\bar{\eta} - | \eta |^2$ are zero. Hence
$\eta$ is the characteristic function of some set. Hence $N = \eta L^2$ is an
invariant subspace contained in $M$.

If $\lambda \in M \ominus N$, then $\lambda$ is orthogonal to $\eta L^2$ and
hence by computing Fourier coefficients $\lambda \bar{\eta}$ is identically
zero. Also
\[ 1 - \eta \perp M \supseteq N \supseteq \{\gamma^n \lambda\} \]
so, computing Fourier coefficients we find $(1 - \bar{\eta}) \lambda$ is
identically zero. Combining these two shows $\lambda$ is the zero function,
hence $M = N$, and the theorem is proved.

\

Clearly, two sets identify the same subspace if and only if their symmetric
difference has zero measure. The Booleian structure of the Borel
$\sigma$-algebra $\mathfrak{B}$ makes the set of the bi-invariant subspaces a
lattice which is isomorphic to $\mathfrak{B}$.

\

We state the characterization of the invariant subspaces of $L^2
(\mathbb{T})$, and sketch its proof.

\begin{theorem}
  The invariant, non-bi-invariant, subspaces of $L^2 (\mathbb{T})$ have the
  form $\psi H^2 (\mathbb{D})$, where $\psi$ is measurable and $| \psi
  (e^{it}) | = 1$ a.e. The function $\psi$ is unique up to a multiplicative,
  unimodular constant.
\end{theorem}

How do we extract $\psi$ from $\mathcal{K}$? For a given invariant subspace
$\mathcal{K}$ such that $S\mathcal{K} \subset \mathcal{K}$, let $\psi \neq 0$
be in $\mathcal{K} \ominus S\mathcal{K} \subseteq \mathcal{K} \ominus S^n
\mathcal{K}$. Then,
\begin{eqnarray*}
  \int_{\mathbb{T}} | \psi (e^{it}) |^2 e^{int} dt & = & < e_n \psi, \psi
  >_{L^2 (\mathbb{T})}\\
  & = & < S^n \psi, \psi >_{L^2 (\mathbb{T})}\\
  & = & 0
\end{eqnarray*}
for $n \geqslant 1$. Similarly $\int_{\mathbb{T}} | \psi (e^{it}) |^2 e^{int}
dt = 0$ for $n \leqslant - 1$, and so $| \psi |$ is a constant, which can be
normalized to $| \psi | = 1$.

\

The reader who is familiar with the spectral theorem can view some of these
results as a special instance of it. The shift is a {\tmem{normal operator}},
$\tau_1^{\ast} \tau_1 = \tau_{- 1} \tau_1 = I = \tau_1 \tau_{- 1}$ (this
implies, more, that $\tau_1$ is an unitary operator on $\ell^2 (\mathbb{Z})$).
Its {\tmem{spectrum}} is $\sigma (\tau_1) =\mathbb{T}$, and the shift can be
identified with the identity map $z \mapsto z$ on $\mathbb{T}$. The measurable
calculus for $\tau_1$ identifies each bounded and Borel measurable $b$ on
$\mathbb{T}$ with the operator $b (\tau_1)$ on $\ell^2 (\mathbb{Z})$; $\sigma
(b (\tau_1)) = \text{ess-range} (b)$, and $\|b\|_{L^{\infty}} = || |b (\tau_1)
| ||$, the operator norm of $b (\tau_1)$. The bi-invariant subspaces of
$\tau_1$ correspond to measurable subsets of the spectrum.

\section{Complex frequencies and the causal case}

\subsection{The Hardy space}

For $\phi \in \ell^2 (\mathbb{N})$, define its $Z$-transform $Z \phi$ to be
\[ Z \phi (z) = \sum_{n = 0}^{\infty} \phi (n) z^n . \]
The series converges to a function holomorphic in the unit disc $\mathbb{D}=
\{z : |z| < 1\}$:
\begin{eqnarray*}
  \left| \sum_{n = M + 1}^N \phi (n) z^n \right|^2 & \leqslant & \sum_{n = M +
  1}^N | \phi (n) |^2 \cdot \sum_{n = M + 1}^N |z|^{2 n}\\
  & \leqslant & \sum_{n = M + 1}^N | \phi (n) |^2 \cdot \frac{|z|^{2 M +
  2}}{1 - |z|^2}
\end{eqnarray*}
which tends to zero uniformly for $|z| \leqslant r < 1$. In holomorphic
control theory the $Z$ transform is usually defined as $Z \phi (z) = \sum_{n =
0}^{\infty} \phi (n) z^{- n}$, and the exterior of the unit disc plays the
role which is in these notes played by the unit disc. What we are doing is
extending the notion of ``frequency'' from $\mathbb{T}$ to $\mathbb{D} \cup
\mathbb{T}$, and the use of the notation $\hat{\phi} (z) = Z \phi (z)$ is
justified.

The old $\hat{\phi} (e^{it})$ can be recovered as the $L^2$.limit of $e^{it}
\mapsto \hat{\phi}  (re^{it})$ as $r \rightarrow 1$,
\begin{eqnarray*}
  \frac{1}{2 \pi}  \int_{\mathbb{T}} | \hat{\phi} (e^{it}) - \hat{\phi}
  (re^{it}) |^2 dt & = \sum_{n = 0}^{\infty} | \phi (n) |^2  (1 - r^{2 n})
  \rightarrow 0 & 
\end{eqnarray*}
as $r \rightarrow 1$.

The \tmverbatim{Hardy space} $H^2 (\mathbb{D})$ is the image of $\ell^2
(\mathbb{N})$ under the $Z$-transform. Alternatively, it can be defined as the
space of the functions $f$ which are holomorphic in $\mathbb{D}$, for which
\[ \|f\|_{H^2}^2 = \sup_{r < 1}  \frac{1}{2 \pi}  \int_{\mathbb{T}} |f
   (re^{it}) |^2 dt = \lim_{r < 1}  \frac{1}{2 \pi}  \int_{\mathbb{T}} |f
   (re^{it}) |^2 dt < \infty . \]
Or, it can be characterized as the space of those $f^{\ast} (e^{it})$ in $L^2
(\mathbb{T})$, $f^{\ast} (e^{it}) = \sum_{n = - \infty}^{+ \infty} \phi (n)
e^{int}$, for which $\phi (n) = 0$ for all negative $n$'s and $\{ \phi (n) \}
\in \ell^2 (\mathbb{Z})$, that is, $\{ \phi (n) \} \in \ell^2 (\mathbb{N})$.
The function $f^{\ast} : \mathbb{T} \rightarrow \mathbb{C}$ is the boundary
function of $f (z) = f (re^{it}) = \sum_{n = 0}^{+ \infty} \phi (n) r^n
e^{int}$, which we identify with $f$, $f = f^{\ast}$.

On the frequency side we have the points of $\mathbb{D}$, and the value of
functions in $H^2$ can be computed at those points, and not just a.e. In fact,
it can be computed in a rather quantitative way.
\begin{eqnarray*}
  f (z) & = & \sum_{n = 0}^{\infty} a_n z^n\\
  & = & < \sum_{n = 0}^{\infty} a_n w^n, \sum_{n = 0}^{\infty} \bar{z}^n w^n
  >_{H^2}\\
  & = & < f (w), \frac{1}{1 - \bar{z} w} >_{H^2}\\
  & = & < f, k_z >_{H^2},
\end{eqnarray*}
where $k (w, z) = k_z (w) = \frac{1}{1 - \bar{z} w}$, $k : \mathbb{D} \times
\mathbb{D} \rightarrow \mathbb{C}$ is the \tmverbatim{reproducing kernel} of
$H^2$.

The theory of {\tmstrong{Hilbert function spaces with a reproducing kernel}}
(RKHS) is old, and it had its inception in work of Bergman and Aronszajn in
the early '40s. Much of what is written in these notes can be proved, or posed
as a problem, for general RKHS's. We will see instances of that in the final
section.

\subsection{The characterization of causal, time invariant systems}

To deal with causal systems, we need $H^{\infty} (\mathbb{D})$, the space of
the bounded analytic functions on the unit disc.

\begin{theorem}
  The causal, time-invariant, linear, $2$-stable systems $T$ are those having
  the form $(T \phi) \widehat{} (z) = b (z)  \hat{\phi} (z)$, with $b$ in
  $H^{\infty}$. Moreover,
  \[ \interleave T \interleave_{\mathcal{B} (H^2)} \assign \sup
     \frac{\|bh\|_{H^2 (\mathbb{T})}}{\|h\|_{H^2 (\mathbb{T})}} =
     \|h\|_{H^{\infty}} . \]
\end{theorem}

Using the maximum principle, it is easy to see that if the transfer function
$b$ is given by the boundary values of a function in $H^{\infty}$, which we
continue to call $b$; then $\sup \frac{\|bh\|_{H^2 (\mathbb{T})}}{\|h\|_{H^2
(\mathbb{T})}} \leqslant \|h\|_{H^{\infty}}$. In the other direction, let $M_b
: H^2 \rightarrow H^2$ be the multiplication operator $f \mapsto bf$, and let
$M_b^{\ast}$ be its adjoint. Then, using the reproducing property of $k_z$,
\begin{eqnarray*}
  M_b^{\ast} k_z (w) & = & < M_b^{\ast} k_z, k_w >\\
  & = & < k_z, M_b k_w >\\
  & = & \overline{< M_b k_w, k_z >}\\
  & = & \overline{< bk_w, k_z >}\\
  & = & \overline{b (z) k_w (z)}\\
  & = & \overline{b (z)} k_z (w),
\end{eqnarray*}
i.e. $M_b^{\ast} k_z = \overline{b (z)} k_z$: the kernel functions are
eigenvectors of the adjoint of the multiplication operator, having the
conjugates of values of $b$ as eigenvalues. This fact holds for general RKHS
and we will encounter it again. We use it now to show the opposite inequality
in the theorem above:
\begin{eqnarray*}
  \sup \frac{\|bh\|_{H^2 (\mathbb{T})}}{\|h\|_{H^2 (\mathbb{T})}} & = & \|M_b
  \|_{\mathcal{B} (H^2)}\\
  & = & \|M_b^{\ast} \|_{\mathcal{B} (H^2)}\\
  & \geqslant & \sup \frac{\|M_b^{\ast} k_z \|_{H^2 (\mathbb{T})}}{\|k_z
  \|_{H^2 (\mathbb{T})}}\\
  & = & \sup | \overline{b (z)} |\\
  & = & \|b\|_{H^{\infty}} .
\end{eqnarray*}
Hidden behind this rather straightforward proof there is a curious fact. There
are $f_{\epsilon}$ in $H^2$ such that
\[ (\|b\|_{H^{\infty}} - \epsilon)  \frac{1}{2 \pi}  \int_{\mathbb{T}}
   |f_{\epsilon} (e^{it}) |^2 dt \leqslant \frac{1}{2 \pi}  \int_{\mathbb{T}}
   |b (e^{it}) f_{\epsilon} (e^{it}) |^2 dt, \]
i.e. $|f_{\epsilon} (e^{it}) |$ is rather concentrated on the set where $|b
(e^{it}) |$ is largest. It is an interesting exercise showing that the
functions $f_{\epsilon}$ can be chosen among kernel functions. (Hint: use the
nonintegrability of $t \mapsto \frac{1}{1 - e^{i t}}$).

The theorem above applies to causal systems having input $\phi$ in $\ell^2
(\mathbb{N})$:
\[ T_b \phi (n) = \sum_{j = 0}^n \check{b}  (n - m) \phi (m), \]
where $\check{b}  (n)$ is the $n^{\tmop{th}}$ coefficient in the series
expansion of $b$ with center at $0$.

The same conclusion applies to $T_b$ having input on the larger space $\ell^2
(\mathbb{Z})$. Passing to the frequency side,
\[ \sup_{\psi \in L^2 (\mathbb{T})}  \frac{\|b \psi \|_{L^2 (\mathbb{T})}}{\|
   \psi \|_{L^2 (\mathbb{T})}} = \sup_{f \in H^2} 
   \frac{\|bf\|_{H^2}}{\|f\|_{H^2}} . \]
In fact, as we have proved, both sides have value $\|b\|_{H^{\infty}} =
\|b\|_{L^{\infty} (\mathbb{T})}$.

We can now give an example of $k \nin \ell^1 (\mathbb{N})$ such that $\phi
\mapsto k \ast \phi$ is bounded on $\ell^2 (\mathbb{N})$. If $k$ were
summable, then $b (z) = \sum_{n = 0}^{\infty} \phi (n) z^n$ would extend to a
function which is continuous on $\bar{\mathbb{D}}$. We only have, then, to
find a bounded holomorphic $b$ which does not admit a continuous extension to
the closed unit disc. As an example, let
\[ b (z) = \exp \left( - \frac{1 + z}{1 - z} \right) . \]
We will see below (and it can be easily verified) that $b$ is {\tmem{inner}}:
bounded and with boundary values of unit modulus $a.e$. The boundary values
are in fact:
\[ b (e^{it}) = \exp \left( \frac{e^{it} + 1}{e^{it} - 1} \right) = \exp (- i
   \cot (t / 2)) \]
which is not continuous at $t = 0$.

This theorem was given a far reaching generalization by von Neumann.

\begin{theorem}
  Let $T$ be a linear contraction on a Hilbert space $H$, $\|Tx\| \leqslant
  \|x\|$, and let $p$ be a complex polynomial. Then,
  \[ \interleave p (T) \interleave \leqslant \|p\|_{H^{\infty}}, \]
  with equality (for any given polynomial $p$) when $H = H^2$ and $T = S$ is
  the shift.
\end{theorem}

This result exemplifies a general trend, of reducing (when possible) questions
concerning a large family of abstract operators to the corresponding question
for a shift-related operator on $H^2$, which works as a {\tmstrong{model}} for
the general theory. A nice reading on these topics is the monograph Nagy and
Fojas (see references).

Observe that the equality $\interleave p (S) \interleave_{\mathcal{B} (H^2)} =
\| p \|_{H^{\infty}}$ holds without restrictions on $p \in H^{\infty}$. In the
general operator theoretic framework this is no longer true.

\subsection{The characterization of the invariant spaces for $\ell^2
(\mathbb{N})$}

A \tmverbatim{inner function} $\Theta$ is a nonconstant function in
$H^{\infty}$ such that $| \Theta (e^{it}) | = 1$ a.e. Such functions play a
preminent role in Hardy theory.

\begin{theorem}
  {\tmstrong{[Beurling]}} The invariant subspaces of $H^2$ have the form
  $\Theta H^2$. The representation is unique up to unimodular constants.
\end{theorem}

Since $H^2 (\mathbb{D})$ is closed in $L^2 (\mathbb{T})$, Beurling's Theorem
easily follows from the characterization of the invariant subspaces for the
shift on $L^2 (\mathbb{T})$. However, the direct approach to the problem is of
interest.

Is is clear that each space having the form $\Theta H^2$ is invariant under
multiplication by $z$. In the opposite direction, we only mention how to find
$\Theta$ if an invariant subspace $\mathcal{K}$ is given. The key point is
showing that $M_z \mathcal{K} \subsetneqq \mathcal{K}$, so we can pick $\Theta
\in \mathcal{K} \ominus M_z \mathcal{K}$ (which will be if necessary
normalized). Let $n \geqslant 0$ be lowest such that $z^n$ divides all $f$ in
$\mathcal{K}$. Then, $n + 1$ is lowest for $M_z \mathcal{K}$, so $M_z
\mathcal{K} \neq \mathcal{K}$.

This simple reasoning, based on the mere existence of a ``order of zero'' for
holomorphic functions, rules out the existence of bi-invariant spaces for the
shift: there are no bi-invariant linear features for signals in positive time.
This is somehow intuitive (the backward shift destroys some of the information
carried by the signal), but it is nonetheless worth mentioning.

The operator $M_{\Theta}$, mapping $H^2$ onto $\Theta H^2$, is an isometry
(but not a unitary operator): $\| \Theta f\|_{H^2} = \|f\|_{H^2}$.

\subsection{The characterization of inner functions}

Since the class of inner functions is the library of ``invariant features'',
it is interesting to have a more concrete characterization for them. There are
two main building blocks we have to consider. The first, generated by Blaschke
products, are determined by the points at which the functions vanish; the
second, the singular inner factors, are determined by the rate at which the
function tends to zero along various radii.

Let $a$ be a point in $\mathbb{D}$. The \tmverbatim{Blaschke factor} $\phi_a
(z) = \frac{|a|}{a}  \frac{a - z}{1 - \bar{a} z}$ maps $\mathbb{D}$,
respectively, $\mathbb{T}$, onto itself, holomorphically and $1 - 1$, hence it
it an inner function. We normalize it so that $\phi_a (a) = 0$ and $\phi_a (0)
= |a| > 0$. Then, the \tmverbatim{finite Blaschke product}
\[ \mathcal{B} (z) = \lambda z^m \Pi_{j = 1}^n  \frac{|a_j |}{a_j}  \frac{a_j
   - z}{1 - \overline{a_j} z}, \]
where $n, m$ are nonnegative integers ($n + m > 0$), $a_1, \ldots, a_n \in
\mathbb{D}$ (repetition being allowed), and $| \lambda | = 1$, is also inner.
It is clear that $\mathcal{B} (z) = 0$ if and only if $z = a_1, \ldots, a_n$
or, if $m > 0$, $z = 0$. In applications to engineering, finite Blaschke
products are especially important, for reasons that will be clear in Section
6. See also the lecture notes of Francis in the reference list.

We can pass to the limit to \tmverbatim{infinite Blaschke products}.

\begin{theorem}
  Let $m$ be a nonnegative integer and $\{a_j \}_{j = 0}^{\infty}$ be a
  sequence in $\mathbb{D}$ (repetition being allowed), and $| \lambda | = 1$.
  Then,
  \[ \mathcal{B} (z) = \lambda z^m \Pi_{j = 1}^{\infty}  \frac{|a_j |}{a_j} 
     \frac{a_j - z}{1 - \overline{a_j} z} \]
  converges to a nonzero holomorphic function in $\mathbb{D}$ if and only if
  the \tmverbatim{Blaschke condition} holds,
  \[ \sum_{j = 1}^{\infty} (1 - |a_j |) < \infty . \]
  Convergence is uniform on compact subsets of $\mathbb{D}$ and $\mathcal{B}
  (z) = 0$ if and only if $z = a_j$ for some $j$, or, if $m > 0$, if $z = 0$.
\end{theorem}

Given a nonconstant, inner function $\Theta$, let $\{a_j \}_{j = 0}^{\infty}$
be the sequence of its zeros $a_j \neq 0$ in $\mathbb{D}$ (repetition being
allowed if the zero has higher order) and let $m \geqslant 0$ be the order of
$\Theta (z)$ at $z = 0$. Then,
\[ \Theta (z) = \lambda \mathcal{B} (z) S (z), \]
where $| \lambda | = 1$, $\mathcal{B} (z) = z^m \Pi_{j = 1}^{\infty} 
\frac{|a_j |}{a_j}  \frac{a_j - z}{1 - \overline{a_j} z}$ is the
\tmverbatim{Blaschke factor} of $\Theta$, normalized to have $\mathcal{B} (0)
> 0$, and $S$ is a inner function with no zero inside $\mathbb{D}$, the
\tmverbatim{singular inner factor} of $\Theta$, $S (0) > 0$.

To have a better understanding of the singular factor, consider the
\tmverbatim{Caley map} $\psi (z) = \frac{1 + z}{1 - z}$, mapping $\mathbb{D}$
one-to-one and onto the right half-plane $\mathbb{C}_+ = \{x + iy : x > 0\}$.
For any $\mu > 0$, the function $S_{0, \mu} (z) = e^{- \mu \psi (z)}$ is then
an inner function, and an $\infty$-one mapping $\mathbb{D}$ onto $\mathbb{D}$
with no zero inside $\mathbb{D}$. \ It tends to zero rapidly as $z = 1 -
\varepsilon$ approaches $1$ along the real axis; $S_{0, \mu}  (1 -
\varepsilon) \sim \exp (- 2 \mu / \varepsilon)$. We might take products of
factors $S_{\alpha, \mu} (z) = S_0, \mu (e^{- i \alpha} z)$ and obtain other
such singular inner functions. We might think of taking infinite products, or
even ``continuous products''. It turns out that such products could well be
``continuous'', but not too much.

\begin{theorem}
  The singular factor has the form:
  \[ S (z) = \exp \left( - \int_{\mathbb{T}} \frac{1 + e^{- it} z}{1 - e^{-
     it} z} d \mu (t) \right), \]
  where $\mu \geqslant 0$ is a Borel measure on $\mathbb{T}$ which is mutually
  singular with respect to arclength measure.
\end{theorem}

When $\mu = \sum_j \mu_j \delta_{\alpha_j}$ is a finite, positive linear
combination of Dirac delta's, then
\[ S (z) = \Pi_j e^{- \mu_j \psi (e^{- i \alpha_j} z)} . \]
At this point we can describe the lattice of (singly) invariant subspaces of
$H^2$. For invariant subspaces generated by Blaschke products the lattice
structure is determined by the lattice of zero sets with the operations $\cap$
and $\cup$. For the subspaces generated by singular functions the lattice is
determined by the lattice of positive singular measures with the operations
$\wedge$ and $\vee$. The full lattice is described by combining these two.

\subsection{Inner/outer factorization}

The multiplication operator $M_{\Theta}$ takes $H^2$ onto the invariant
subspace $\Theta H^2$. It turns out that all multiplication operators we have
seen in the analysis of causal systems admit a canonical factorization through
an operator of this sort. Actually, it is convenient to look at things in more
generality.

A function $u$ in $H^1$ is \tmverbatim{outer} if
\[ u (z) = \exp \left( \frac{1}{2 \pi}  \int_{\mathbb{T}} \frac{1 + e^{- it}
   z}{1 - e^{- it} z} k (e^{it}) dt \right), \]
for some real valued, integrable $k$ on $\mathbb{T}$. The function $k$ can be
easily recovered from $u$:
\[ k (e^{it}) = \log |u (e^{it}) |, a.e. \]
We have chosen a normalization for which $u (0) > 0$.

\begin{theorem}
  \label{pippo}Let $b$ be in $H^1$. Then, there are a unique outer function
  $u$ and inner function $\Theta$ such that
  \[ b = u \Theta . \]
  Moreover, $\|b\|_{H^p} = \|u\|_{H^p}$ for $p = 1 \leqslant p \leqslant
  \infty$.
\end{theorem}

Outer functions $u \in H^{\infty} (\mathbb{D})$ can be characterized as those
which are invertible in the weak sense that $u H^2 (\mathbb{D})$ is dense in
$H^2 (\mathbb{D})$. In fact, more can be said.

\begin{theorem}
  Let $f$ be in $H^2$ and let $[f]$ be the smallest invariant subspace of
  $H^2$ containing $f$. Then, with $\Theta u$ as in the inner/outer
  factorization of $f$, we have
  \[ [f] = \Theta H^2 . \]
\end{theorem}

\qquad Hence if $f$ is outer then $[f] = H^2$ and in particular $1 \in [f]$.
Thus $f$ is invertible in $H^2$ in the weak sense that there is a sequence $\{
g_n \} \subset H^2$ such that $g_n f \rightarrow 1$ in the norm of $H^2$.
However $1 / f$ need not be in $H^2$; for instance $f (z) = 1 - z$ is outer
(as is most easily seen by computing $[1 - z]^{\perp}$, i.e. showing that $H^2
(\mathbb{D}) \ominus (1 - z) H^2 (\mathbb{D}) = 0$). Inner functions are not
invertible in $H^{\infty}$; further, if $\Theta$ is inner then $[\Theta] =
\Theta H^2 \subsetneqq H^2$ and thus $\Theta$ does not even have an inverse in
a weak sense we just saw.

\qquad Thus if $b$ has the inner/outer factorization $b = \Theta u$ then we
can write the operator $M_b$ as a product of two commuting operators; the
isometric map $M_{\Theta}$ which imposes"features" on the signal, and $M_u$
which is a (roughly) invertible operator on the space of functions with
specified features.

Another consequence of the inner/outer factorization is the following.

\begin{lemma}
  \label{pluto}For $h \in H^1$ we have
  \[ \|h\|_{H^1} = \inf \{\|f\|_{H^2} \|g\|_{H^2} : h = fg\} . \]
\end{lemma}

The $\leqslant$ direction is just Cauchy-Schwarz. In the other direction, we
can write $h = u \Theta$ with $u$ outer, then zero free in $\mathbb{D}$: $h =
(u^{1 / 2})  (u^{1 / 2} \Theta) = fg$, with $\|h\|_{H^1} = \|f\|_{H^2}
\|g\|_{H^2}$.

\section{Approximating noncausal systems by causal ones: Hankel operators and
Nehari theory}

Given a function $\phi \in L^{\infty} (\mathbb{T})$, here identified with the
invariant operator $\psi \mapsto M_{\phi} \psi = \phi \psi$ on $L^2
(\mathbb{T})$, what is the best approximation of $M_{\phi}$ by causal
operators $M_b$ with $b \in H^{\infty}$? Namely, we look for
\[ \inf_{b \in H^{\infty}} \sup_{f \in H^2}  \frac{\| \phi f -
   bf\|_{H^2}}{\|f\|_{H^2}} = \inf_{b \in H^{\infty}} \| \phi -
   b\|_{L^{\infty}} = \tmop{dist} (\phi, H^{\infty}) . \]
Indeed, one would also like to know if a minimizing $b$ exists (yes), if it is
unique (sometimes, in many relevant cases), if there is a way to construct it
(again, yes in many cases of interest).

In the passage from first to second member the $\leqslant$ direction is
obvious. For the opposite direction, note that the $L^{\infty}$ norm of $\phi
- b$ requires testing on $L^2$ functions, while on the left we only test on
$H^2$ functions. We use the shift invariance of the $L^2 (\mathbb{T})$ norm.
For $\epsilon > 0$ let $\psi \in L^2 (\mathbb{T})$ be such that $\| \psi
\|_{L^2} = 1$ and $\| \phi \psi \|_{L_2} \geqslant \| \phi \|_{L^{\infty}} -
\epsilon$. Find $N$ such that for $|z| = 1$, $\psi_N (z) = \sum_{n = -
N}^{\infty} \hat{\psi} (n) z^n$ satisfies $\| \psi - \psi_N \|_{L^2} <
\epsilon$. Then,
\begin{eqnarray*}
  \| \phi \|_{L^{\infty}} - \epsilon & \leqslant & \| \phi \psi \|_{L_2}\\
  & \leqslant & \| \phi \psi_N \|_{L_2} + \| \phi (\psi - \psi_N)\|_{L_2}\\
  & \leqslant & \left( \frac{1}{2 \pi}  \int_{\mathbb{T}} | \phi (e^{it})
  \psi_N (e^{it}) |dt \right)^{1 / 2} + \| \phi \|_{L^{\infty}} \cdot
  \epsilon\\
  & = & \left( \frac{1}{2 \pi}  \int_{\mathbb{T}} | \phi (e^{it}) e^{iNt}
  \psi_N (e^{it}) |dt \right)^{1 / 2} + \| \phi \|_{L^{\infty}} \cdot
  \epsilon\\
  & = & \left( \frac{1}{2 \pi}  \int_{\mathbb{T}} | \phi (e^{it}) f (e^{it})
  |dt \right)^{1 / 2} + \| \phi \|_{L^{\infty}} \cdot \epsilon
\end{eqnarray*}
where $f (z) = z^N \psi_N (z)$ is holomorphic and $1 \geqslant \|f\|_{H^2} =
\| \psi_N \|_{L^2}$. Thus,
\[ \frac{\| \phi f\|_{L^2}}{\|f\|_{L^2}} \geqslant \| \phi \|_{L^{\infty}}  (1
   - \epsilon) - \epsilon, \]
and the $\geqslant$ direction in the equality is proved. A shorter proof can be
derived using Toeplitz operators.

\

\subsection{Hankel forms and Hankel operators}

The approximation problem just described, finding $b$, the optimal
$H^{\infty}$ approximation to $\phi$, can be stated in the language of Hankel
operators and Nehari's theorem characterizing the norm of Hankel gives
information about $b$. We begin with some definitions.

The \tmtexttt{Hankel matrix operator} $\Gamma_{\alpha}$ induced by a complex
valued sequence $\alpha = \{\alpha_n \}_{n = 0}^{\infty}$ is defined on
sequences $a = \{a_n \}_{n = 0}^{\infty}$ (in $\ell_c (\mathbb{N})$, to start
with) by
\[ (\Gamma_{\alpha} a) (m) = \sum_{n = 0}^{\infty} \alpha_{m + n} a_n, \]
or
\[ < \Gamma_{\alpha} a, b >_{\ell^2} = \sum_{m.n \geqslant 0} \alpha_{m + n}
   a_n \overline{b_m} \]
A famous example of an Hankel matrix is Hilbert's matrix $[(i + j + 1)^{-
1}]_{i, j = 0}^{\infty}$.

We have already seen how useful it is to pass to the frequency side by the
$Z$-transform. Let $P_+$ be the orthogonal projection of $L^2 (\mathbb{T})$
onto $H^2$ and for any $g \in L^2 (\mathbb{T})$ write $g_+ = P_+ g$ and $g_- =
g - g_+$. Hence $g_-$ is the projection of $g$ onto $L^2 \ominus H^2$ and the
$g_-$ obtained this way are exactly the functions $\overline{zj}$ for $j \in
H^2$. For $\phi \in L^2 (\mathbb{T})$ we define the \tmtexttt{Hankel bilinear
form} $B_{\phi}$ associated to $\phi$, a bilinear map $H^2 \times H^2
\rightarrow \mathbb{C}$ and define the \tmtexttt{Hankel operator} with symbol
$\phi$, $H_{\phi}$, to be the linear map of $H^2$ to $L^2 \ominus H^2$ by
\[ B_{\phi}  (f, g) \assign \langle fg, \bar{z}  \bar{\phi} \rangle_{L^2} = :
   \langle H_{\phi} f, \overline{zg} \rangle_{L^2} . \]
In particular $H_{\phi} f = (\phi f)_-$.

The relation between Hankel forms and Hankel matrices is the following:

\begin{align*}
  B_{\phi}  (f, g) & = \langle \sum_{n = 0}^{\infty} \hat{f} (n) z^n  \sum_{m
  = 0}^{\infty} \hat{g} (m) z^m, \sum_{k = - \infty}^{\infty}
  \overline{\hat{\phi} (k)} z^{- k - 1} \rangle_{L^2}\\
  & = \sum_{k \leqslant - 1} \overline{\hat{\phi} (k)} \sum_{m + n = - k -
  1}^{\infty} \hat{f} (n)  \hat{g} (m)\\
  & = \sum_{m = 0}^{\infty} \hat{g} (m)  \sum_{n = 0}^{\infty}
  \overline{\hat{\phi}  (- m - n - 1)} \hat{f} (n)\\
  & = \langle \Gamma_{\alpha}  \hat{f}, \overline{\hat{g}} \rangle_{\ell^2},
\end{align*}

where $\alpha (j) = \overline{\hat{\phi}  (- j - 1)}$. From these formal
calculations it is evident that
\[ [B_{\phi}] \assign \sup_{f, g \in H^2}  \frac{|B_{\phi} (f, g)
   |}{\|f\|_{H^2} \|g\|_{H^2}} = \| H_{\phi} \|_{\mathtt{operator}} =
   \interleave \Gamma_{\alpha} \interleave_{\mathcal{B} (\ell^2)} . \]

If $\gamma$ is bounded then
\[ | B_{\gamma} (f, g) | = | \langle fg, \bar{z}  \bar{\gamma} \rangle_{L^2} |
   = \left| \frac{1}{2 \pi}  \int_{\mathbb{T}} f (e^{it}) g (e^{it}) e^{it}
   \gamma (e^{it}) dt \right| \leqslant \| \gamma \|_{L^{\infty}} \|f\|_{H^2}
   \|g\|_{H^2}, \]
and hence $[B_{\gamma}] \leq \| \gamma \|_{L^{\infty}}$.Also clearly for any
$b \in H^2$ $B_{\phi} = B_{\phi - b}$. Combining these facts we have
\[ [B_{\phi}] \leqslant \inf \{\| \phi - h\|_{L^{\infty}} : h \in H^2 \} =
   \mathrm{dist} (\phi, H^{\infty}) . \]
Given $\phi \in L^2$ let $b \in H^2$ be that function, if there is one, such
that $\| \phi - b\|_{L^{\infty}} = \mathrm{dist} (\phi, H^{\infty})$. If
$\phi$ is bounded then $b$ is in $H^{\infty}$ and is the function we discussed
earlier, the best approximation to $\phi$ in the $L^{\infty}$ norm. To
complete the story we show the opposite inequality, and will then know that
the norm of the Hankel operator, or of the Hankel form, equals the distance of
the symbol from $H^{\infty}$. That result is Nehari's theorem.

\begin{theorem}
  Given $\phi \in L^2$
  \[ [B_{\phi}] = \| H_{\phi} \|_{\mathtt{operator}} = \interleave
     \Gamma_{\alpha} \interleave_{\mathcal{B} (\ell^2)} = \mathrm{dist} (\phi,
     H^{\infty}) . \]
\end{theorem}

The previous discussion shows that the expression on the right is larger. To
finish we must show that there is a holomorphic function $b$ so that $\| \phi
- b \|_{\infty} = [B_{\phi}]$. Starting with the formula $B_{\phi}  (f, g)
\assign \langle fg, \bar{z}  \bar{\phi} \rangle_{L^2}$ and taking note of
Lemma 1 which shows that $fg$ is a generic element of $H^1$ we see that
$[B_{\phi}]$ is equal to the norm of the functional $h \rightarrow \langle fg,
\bar{z}  \bar{\phi} \rangle_{L^2}$ acting on $H^1$. By the Hahn-Banach theorem
that functional extends in a norm preserving way to a functional on $L^1$.
That functional on $L^1$ will be of the form $k \rightarrow \langle k, j
\rangle_{L^2}$ for a bounded $j$ with $\| j \|_{\infty} = [B_{\phi}]$ and $j$
will satisfy
\[ \langle h, \bar{z}  \bar{\phi} \rangle_{L^2} = \langle h, j \rangle_{L^2}
   \text{ } \forall h \in H^1 . \]
In particular $j$ and $\bar{z}  \bar{\phi}$ have the same nonnegative Fourier
coefficients and thus $j_= = (\bar{z}  \bar{\phi})_+$.We now want to find $b$
so that $\| \phi - b \|_{\infty} = \| j \|_{\infty}$. We have
\[ \overline{z \phi} = (\overline{z \phi})_+ + (\overline{z \phi})_- = j_+ +
   (\overline{z \phi})_- = j - j_- + (\overline{z \phi}) . \]
Rearranging gives $\overline{z \phi} - (- j_- + (\overline{z \phi})) = j$.
From that one quickly shows there is a holomorphic $b$ so that $\phi - b = z
\bar{j}$ and that is enough to give what we want, because $\| \phi - b
\|_{L^{\infty}} = \| j \|_{L^{\infty}} = [B_{\phi}] \leqslant \| \phi - b
\|_{L^{\infty}}$.

On Hankel operators, for the mathematical side a good starting point is
Peller's survey; their use in control theory is in Francis' lecture notes.

\subsection{Detour: Toeplitz operators}

For $\psi \in L^2 (\mathbb{T})$ given, the Toeplitz operator $T_{\psi}$ with
symbol $\psi$ is defined for $f \in H^2$ by $T_{\psi} f = M_{\psi} f -
H_{\psi} f = P_+  (\psi f)$, where $P_+ : L^2 \rightarrow H^2$ is orthogonal
projection. The Toeplitz operator coincides with the multiplication operator
$M_{\psi}$ if $\psi \in H^2$ is holomorphic. The adjoint of $T_{\psi}$ is
$T_{\psi}^{\ast} = T_{\bar{\psi}}$.

In signal theory Toeplitz operators naturally appear in connection with an
alternative definition of \ on $\ell^2 (\mathbb{N})$. Recall that $\tau_1 \phi
(n) = \phi (n - 1)$ defines the shift on $\ell^2 (\mathbb{N})$. Its adjoint,
the \tmverbatim{backward shift}, is the operator $\tau_1^{\ast} \phi (n) =
\phi (n + 1)$, $\tau_1^{\ast} : \ell^2 (\mathbb{N}) \rightarrow \ell^2
(\mathbb{N})$. It is readily verified that $\tau_1^{\ast} \tau_1 \phi = \phi$
and that $\tau_1 \tau_1^{\ast} \phi = \phi - \phi (0) \delta_0$. A linear
system $T$ on $\ell^2 (\mathbb{N})$ is called time invariant if $\tau_1^{\ast}
T \tau_1 = T$: if we shift the input forward, feed it to $T$, then shift
backward, we have the same as just applying $T$.

The rationale for this new definition of invariant system for signals in
positive time is that the previous definition assumed, in order to be
verified, that all the past values of the signal have been stored and are
accessible, a requirement which is not practical.

We now see how invariant systems lead to Toeplitz operators. Passing to the
frequency side with $h (z) = \sum_{j = 0}^{\infty} a_j z^j$, a (linear) system
$T$ on $\ell^2 (\mathbb{N})$, represented by a matrix $[F_{ij}]_{i, j =
0}^{\infty}$ ($F_{ij} = < T (z^j), z^i >_{H^2}$ are the matrix elements of $T$
with respect to the basis $\{z^n \}_{n = 0}^{\infty}$ of $H^2$), is invariant
if
\begin{eqnarray*}
  \sum_{i = 0}^{\infty} \sum_{j = 0}^{\infty} F_{ij} a_j z^i & = & (T \phi)
  \widehat{} (z)\\
  & = & (\tau_1^{\ast} T \tau_1) \widehat{} \phi (z)\\
  & = & \bar{z}  ((T \tau_1) \widehat{} \phi (z) - (T \tau_1) \widehat{} \phi
  (0))\\
  & = & \bar{z}  \left( \sum_{i = 0}^{\infty} \sum_{j = 1}^{\infty} F_{ij}
  a_{j - 1} z^i - \sum_{j = 1}^{\infty} F_{0 j} a_{j - 1} \right)\\
  & = & \sum_{i = 1}^{\infty} \sum_{j = 1}^{\infty} F_{ij} a_{j - 1} z^{i -
  1}\\
  & = & \sum_{i = 0}^{\infty} \sum_{j = 0}^{\infty} F_{i + 1, j + 1} a_j z^i,
\end{eqnarray*}
i.e. $F_{i + 1, j + 1} = F_{i, j}$: $T$ is represented, w.r.t. the basis
$\{z^n \}_{n = 0}^{\infty}$, by a \tmverbatim{Toeplitz matrix} $F_{i, j} =
f_{i - j}$. Recall that in a Hankel matrix the $i, j$ entry is a function of
$i + j$.

Inserting this back in the expression for $T$ in frequency space,
\begin{eqnarray*}
  (T \phi) \widehat{} (z) & = & \sum_{i = 0}^{\infty} \sum_{j = 0}^{\infty}
  f_{i - j} a_j z^i\\
  & = & P_+  (b (z) f (z)),
\end{eqnarray*}
where $\psi (z) = \sum_{n = - \infty}^{+ \infty} f_n z^n$.

When $\phi$ is holomorphic, the matrix $[f_{i - j}]$ is lower triangular.

As with Hankel operators, it is clear that $\interleave T_{\phi}
\interleave_{\mathcal{B} (H^2)} \leqslant \| \phi \|_{L^{\infty}}$: $\|P_+
(\phi f)\|_{H^2} \leqslant \| \phi f\|_{L^2} \leqslant \|f\|_{H^2}$. Contrary
to the Hankel case, there is no way to improve this estimate:
\[ \interleave T_{\phi} \interleave_{\mathcal{B} (H^2)} = \| \phi
   \|_{L^{\infty}} . \]
Let $k_a (z) = \frac{1}{1 - \bar{a} z}$ be the reproducing kernel at $a$:
\begin{eqnarray*}
  | \langle T_{\phi} k_a, k_a \rangle | & = & | \langle P_+ (\phi k_a), k_a
  \rangle |\\
  & = & | \langle \phi k_a, k_a \rangle |\\
  & = & \left| \frac{1}{2 \pi}  \int_{- \pi}^{\pi} \phi (e^{it}) |k_a
  (e^{it}) |^2 dt \right|\\
  & = & \frac{1}{1 - |a|^2}  \left| \frac{1}{2 \pi}  \int_{- \pi}^{\pi} \phi
  (e^{it}) \frac{1 - |a|^2}{|1 - \bar{a} z|^2} dt \right|\\
  & = & \|k_a \|_{H^2}^2  |P \phi (a) |,
\end{eqnarray*}
where $P \phi$ is the Poisson integral of $\phi$ at $a$, because $P (a,
e^{it}) = \frac{1}{2 \pi}  \frac{1 - |a|^2}{|1 - \bar{a} z|^2}$ is the Poisson
kernel in the unit disc. Hence,
\begin{eqnarray*}
  \interleave T_{\phi} \interleave_{\mathcal{B} (H^2)} & \geqslant & \sup_{a
  \in \mathbb{D}}  \frac{| \langle T_{\phi} k_a, k_a \rangle |}{\|k_a
  \|_{H^2}^2}\\
  & = & \|P \phi \|_{L^{\infty} (\mathbb{D})}\\
  & = & \| \phi \|_{L^{\infty} (\mathbb{T})} .
\end{eqnarray*}

\subsection{$H^1$ and $\tmop{BMO}$}

We will not touch here \tmverbatim{Nehari's problem}; that is, how to find the
best approximant of $\phi$ in $H^{\infty}$. Even the estimate we have found,
however, are of little use unless we have tools for estimating
$\|b\|_{(H^1)^{\ast}}$. Contrary to a first, naif guess, the dual of $H^1$
contains, but is larger, than $H^{\infty}$.

Shortly after Nehari's article on Hankel forms, Fritz John introduced, in
connection to problems in elasticity theory, the space $\tmop{BMO}$ of
functions with Bounded Mean Oscillations, which he further studied together
with John Nirenberg. Restricted to functions on $\mathbb{T}$, the definition
is as follows. For each arc $I \subset \mathbb{T}$, denote by $\phi_I =
\frac{1}{|I|}  \int_I \phi (e^{it}) dt$ be the average of $\phi$ over $I$. The
mean oscillation of $\phi$ over $I$ is $\frac{1}{|I|}  \int [\phi (e^{it}) -
\phi_I] dt, \tmop{andthe}$BMO \ norm of $\phi$ is
\[ \| \phi \|_{L^{\infty}} + \sup_I  \frac{1}{|I|}  \int | \phi (e^{it}) -
   \phi_I | dt. \]
In 1971 C. Fefferman made the surprising discovery that $(H^1)^{\ast} =
\tmop{BMOA}$, the space of the $\tmop{BMO}$ functions which extend
holomorphically to the unit disc. Duality is with respect to the $H^2$ inner
product. It is not difficult to see that this result implies that if $\phi$ is
bounded, then $H \phi$, its Hilbert transform, belongs to $\tmop{BMO}$.

On his way to the proof, Fefferman proved that the $\tmop{BMO}$ norm of a
function can be characterized in terms of \tmverbatim{Carleson measures}. Let
$\mu \geqslant 0$ be a Borel measure on $\mathbb{D}$. We say that it is a
Carleson measure for $H^2$ if there is a positive constant $[\mu]_{\tmop{CM}}$
such that
\[ \int_{\mathbb{D}} |f|^2 d \mu \leqslant [\mu]_{\tmop{CM}} \|f\|_{H^2}^2 .
\]
The concept itself had been introduced by Carleson in connection to the
problem of interpolating functions in $H^{\infty}$. Fefferman showed that $b
\in \tmop{BMOA}$ if and only if $d \mu_b (z) = (1 - |z|^2) |b' (z) |^2 dxdy$
is a Carleson measure.

The appearence of such measures is easily explained. A equivalent norm for
$H^2$ is
\[ [f]_{H^2}^2 = |f (0) |^2 + \int_{\mathbb{D}} (1 - |z|^2) |f' (z) |^2 dxdy.
\]
If $d \mu_b$ is a Carleson measure for $H^2$, then (assuming momentarily that
$b (0) = 0$ and using the equivalent norm to define the inner product),
\begin{eqnarray*}
  | \langle fg, b \rangle_{H^2} | & = & \left| \int_{\mathbb{D}} (fg)'
  \overline{b'} (1 - |z|^2) dxdy \right|\\
  & \leqslant & \left| \int_{\mathbb{D}} fg' \overline{b'} (1 - |z|^2) dxdy
  \right| + \left| \int_{\mathbb{D}} gf' \overline{b'} (1 - |z|^2) dxdy
  \right| \leqslant\\
  & \leqslant & \left| \int_{\mathbb{D}} |g' |^2 (1 - |z|^2) |b' |^2 dxdy
  \right|^{1 / 2} \|f\|_{H^2} + \left| \int_{\mathbb{D}} |f' |^2 (1 - |z|^2)
  |b' |^2 dxdy \right|^{1 / 2} \|g\|_{H^2}\\
  & = & \left| \int_{\mathbb{D}} |g' |^2 d \mu_b \right|^{1 / 2} \|f\|_{H^2}
  + \left| \int_{\mathbb{D}} |f' |^2 d \mu_b \right|^{1 / 2} \|g\|_{H^2}\\
  & \leqslant & 2 [\mu_b]_{\tmop{CM}} \|f\|_{H^2} \|g\|_{H^2} .
\end{eqnarray*}
Recalling Section 5.1, this shows that if $\mu_b$ is Carleson, then the Hankel
form $B_b$, hence the Hankel operator $H_b$, is bounded. By Nehari's theorem,
$b \in (H^1)^{\ast}$. The delicate point is proving the opposite implication.

The short and dense monograph of Sarason well explains the connections between
Hankel operators, basic questions of operator theory, and harmonic analysis.

We summarize part of what we have seen in a diagram:
\[ \tmop{Mult} (H^2) = H^{\infty} \hookrightarrow \tmop{BMOA} = (H^1)^{\ast}
   \hookrightarrow H^2 \hookrightarrow H^1 = H^2 \cdot H^2 . \]
We see here, as it often happens, that analysis on a function Hilbert space
requires introducing a number of other Banach function spaces.

\section{Systems and feedback}

A typical device (a \tmverbatim{plant}) can be modeled by a linear, time
invariant, causal, stable operator $P$, which acts in frequency as $M_b$, with
$b \in H^{\infty}$, and which we assume to be free of feedback loops.
Generally the output $Pa (n)$ only depends on finitely many values $a_{n - m +
1}, \ldots, a_n$ of the input (which have to be stored), and it is easy to
verify that this holds if and only if $b$ is a polynomial of degree $m$. This
property is sometimes expressed saying that transient inputs produce transient
outputs, and it is clear that it suffices to verify this for the unit impulse
$\delta_0$.

A {\tmem{feedback system}} is one in which the output of $P$ is ``fed back''
into $P$, possibly after having been processed by a different plant $C$. For
instance:

\

\resizebox{10.5947297651843cm}{4.49578577987669cm}{\includegraphics{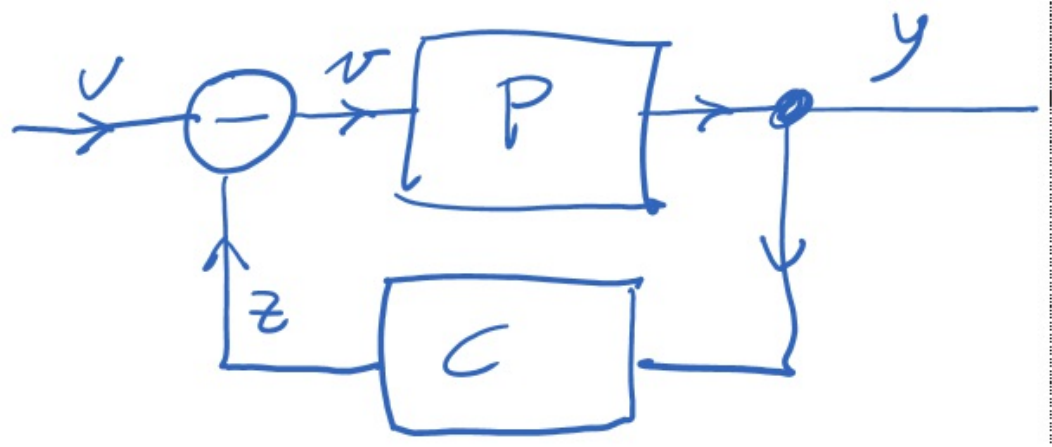}}

We use the same symbols for signals and plants, and their $Z$-transforms and
transfer functions; the letter $n$ stands for time and $\omega$ for frequency.
In a real situation, the output $y (n)$ can not immediately affect the input
$u (n)$ at time $n$. In order to have this, $C (\omega)$ must include a delay
by at least a time unit; i.e. the polynomial $C (0) = 0$.

The system represented by the diagram is:
\[ \left\{ \begin{array}{l}
     y = Pv\\
     v = u - z\\
     z = Cy
   \end{array} \right. \]
Overall, $y = P (u - Cy)$, i.e. $y = \frac{P}{1 + PC} u$. Observe that the
rational function $\frac{P}{1 + PC}$ is not a polynomial, hence the system
with feedback gives a persistent signal as output if the input is the unit
impulse (the feedback produces an ``echo'').

This easy example shows how nontrivial conclusions can be drawn by elementary
algebra in frequency space. Hardy space theory leads system theory much
further. We give here just one example, giving us the opportunity of
mentioning Pick theory, a topic of current research.

\subsection{The model matching problem and the Pick property}

Let $T$ be an ideal plant we want to best approximate by a cascade $UCV$,
where $U$ and $V$ are given plants, and $C$ is a plant we can design. That
is, we want to find $C$ which minimizes
\[ \|T - UCV\|_{H^{\infty}} . \]
This is the {\tmem{\tmverbatim{Model Matching Problem}}} with data $T, U, V$.

\resizebox{11.7448183130001cm}{4.67004460186278cm}{\includegraphics{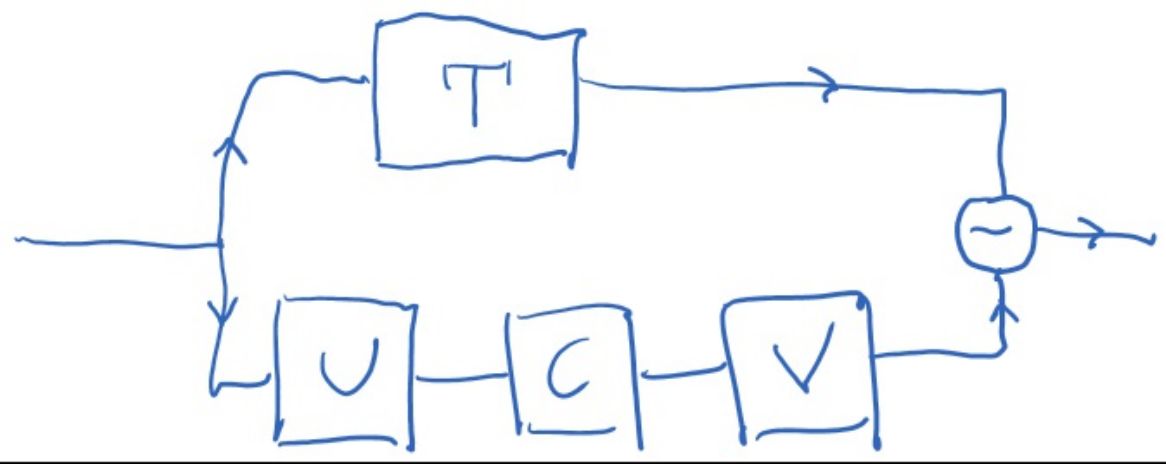}}

Consider the inner/outer factorization $UV = A_i A_o$. Since $U$ and $V$ are
rational (we allow feedbacks), $A_i$ is a finite Blaschke product, with zeros
$\lambda_1, \ldots, \lambda_n$ in $\mathbb{D}$. We can then write $H \assign T
- UCV = T - A_i F$, with $F = A_o C$. Since $H (\lambda_j) = \mu_j \backassign
T (\lambda j)$, we have that (assuming a minimizer exists):
\[ \min \{\|T - A_i F\|_{H^{\infty}} : F \in H^{\infty} \} = \min
   \{\|H\|_{H^{\infty}} : H (\lambda_j) = \mu_j \} . \]
In fact, if $H$ is a minimizer for the right hand side, then the equation $T -
A_i F = H$ has the solution $F = A_i^{- 1}  (T - H)$, which is well defined in
$H^{\infty}$ because $T - H$ vanishes at $\lambda_1, \ldots, \lambda_n$. Since
$A_o$ is outer, we can then reconstruct $C$ from $F$.

Finding a function $H$ of minimal norm satisfying the interpolation costraints
$H (\lambda_j) = \mu_j$ is the {\tmem{\tmverbatim{Pick problem}}} with data
$\{\lambda_1, \ldots, \lambda_n \}$ and $\{\mu_1, \ldots, \mu_n \}$.

Suppose that the minimal norm of $H$ is not larger that $R$. We have sequences
$\{\lambda_1, \ldots, \lambda_n \}$ and $\{\mu_1 / R, \ldots, \mu_n / R\}$ in
$\mathbb{D}$ and we have an interpolating $H / R$ of norm at most one. A
necessary condition for this to hold is the \tmverbatim{Pick property}. For
any choice of complex $a_1, \ldots, a_n$, denoting by $k_{\lambda_j}$ the
reproducing function at $\lambda_j$ and using $M_H^{\ast} (k_{\lambda_j}) =
\overline{H (\lambda_j)} k_{\lambda_j}$:
\begin{eqnarray*}
  0 & \leqslant & \left| \sum_{j = 1}^n a_j k_{\lambda_j} \right|^2 - \left|
  \frac{1}{R} M_H^{\ast} \left( \sum_{j = 1}^n a_j k_{\lambda_j} \right)
  \right|^2\\
  & = & \left| \sum_{j = 1}^n a_j k_{\lambda_j} \right|^2 - \left|
  \frac{1}{R}  \sum_{j = 1}^n a_j \overline{H (\lambda_j)} k_{\lambda_j}
  \right|^2\\
  & = & \left| \sum_{j = 1}^n a_j k_{\lambda_j} \right|^2 - \left| \sum_{j =
  1}^n a_j \overline{\mu_j} / Rk_{\lambda_j} \right|^2\\
  & = & \sum_{i, j = 1}^n a_i \overline{a_j} < k_{\lambda_i}, k_{\lambda_j}
  >_{H^2} \left( 1 - \frac{\overline{\mu_i} \mu_j}{R^2} \right)\\
  & = & \sum_{i, j = 1}^n a_i \overline{a_j} k_{\lambda_i} (\lambda_j) 
  \left( 1 - \frac{\overline{\mu_i} \mu_j}{R^2} \right)
\end{eqnarray*}
That is, the \tmverbatim{Pick matrix} $\left[ k_{\lambda_i} (\lambda_j) \left(
1 - \frac{\overline{\mu_i} \mu_j}{R^2} \right) \right]_{i, j = 1}^n$ is
positive semidefinite.

Pick's Theorem says that the converse is true.

\begin{theorem}
  Given points $\lambda_1, \ldots, \lambda_n$ in the unit disc and values
  $\mu_1, \ldots, \mu_n$ in the unit disc, there exists a function $H$ in
  $H^2$ having norm at most one interpolating them, $H (\lambda_j) = \mu_j$,
  if and only if the matrix
  \[ [k_{\lambda_i} (\lambda_j) (1 - \overline{\mu_i} \mu_j)]_{i, j = 1}^n \]
  is positive definite.
  
  Moreover, the interpolating function $H$ of minimal norm is a rational
  function.
\end{theorem}

Pick's Theorem holds, with natural modifications, for infinite sequences of
points and values.

Extensions and applications of Pick theory are one of the most active areas of
current research at the frontier between operator theory and function spaces.

\section{Beyond the Hardy space; RKHS}

Up to this point this article could be seen as a bus tour of an interesting
city. The bus goes from place to place, the tour guide offers enthusiastic
description and commentary, and at a few of the places the passengers have a
chance to get off the bus and look in detail at some of the sights. That tour
is over now and what comes next can be seen as the airplane ride home. We fly
over a landscape and a voice on the speaker points out some interesting
features below; just a quick glance at them, perhaps enough to whet the
appetite.

The Hardy space lives in the intersection of several powerful mathematical
technologies. Hardy space functions are holomorphic functions in the disk and
can be analyzed using tools from function theory. The boundary values of Hardy
space functions are in the Lebesgue space $L^2$ of the circle and hence the
machinery of Fourier analysis can be used can be used to study them. In fact
they form a closed subspace of $L^2$ and hence there is an associated
projection operator and that lets questions about Hardy space functions be be
formulated and studied in the language of linear operators on Hilbert space.
We have seen bits of all of these approaches.

A point of view we are emphasizing here is that the Hardy space is a Hilbert
space with reproducing kernel, RKHS. That is, it is a Hilbert space whose
elements are functions on a set $X$ (in this case $X = \mathbb{D)}$, the
evaluation of the functions at points $x \in X$ are continuous linear
functionals, and hence each of those evaluations is given by taking the inner
product of $f$ with some distinguished element $k_x$ in the space; $f (x) =
\langle f.k_x \rangle$. The $k_x$ are the reproducing kernels and, in some
sense, the collection of them, $\{ k_x \}_{x \in X}$, plays the role in this
theory that an orthonormal basis plays for finite dimensional inner product
spaces.

In the next three sections we take a look, from great height, at three other
examples of RKHS. The first is the Paley Wiener space, a space somewhat
similar to the Hardy space (but the Hardy space of the half plane rather than
the disk) that is of great interest in the theory of sampling and
reconstructing band limited signals such as speech and music. The second
example is the Dirichlet space. It is a variation of the Hardy space with some
similarities and some differences, and it is dear to the authors. The third
example is the dyadic Dirichlet space. That space is a simplified model of the
Dirichlet space, useful in analyzing the Dirichlet space. It is also a space
which makes explicit the parameter space for the "phase space analysis" of
signals, of which wavelets are the most prominent example.

\subsection{Paley Wiener Space}

The Paley-Wiener space, $PW$, is the subspace of $L^2 (\mathbb{R)}$ of all
functions $f$ whose Fourier transform $\hat{f}$ supported on the interval $[-
\pi, \pi]$, The space is often used in signal analysis; $f \in PW$ is a
signal, $f (t)$ is its value at time $t \in \mathbb{R}$ is its value at time
$t$ and $\hat{f}$, its Fourier transform is the frequency space representation
of the signal. The fact that $\hat{f}$ is supported in $[- \pi, \pi]$ is is a
statement that the signal contains no frequencies outside this range, the
signal is "band limited". The norm of $f$ in $PW$, which is the same as the
norm of $f$ in $L^2$ and (with our normalization) the same as the norm of
$\hat{f}$ in $L^2 (- \pi, \pi)$, is the energy of the signal. In short $PW$ is
a space of finite energy band limited signals. This can be compared with the
Hardy space of the upper half plane; the boundary values of those functions
are exactly the functions $f$ with $\hat{f} \in L^2 (0, \infty)$.

(The same space of functions can also be defined by restricting to the real
axis a certain class of entire functions defined by their growth at infinity.
The equivalence between the two definitions uses the fundamental ideas
developed by Paley and Wiener in the 1930's relating the smoothness of
functions and the decay of their Fourier transforms.)

To see that $PW$ is an RKHS we want to know that the evaluations of points of
$\mathbb{R}$ are continuous functionals. Consider first evaluation at $t = 0$.
We now describe the picture from our very high altitude. The value of $f$ at
$0$ is gotten by using the bilinear pairing $(f, g) \rightarrow \int f
\bar{g}$ to pair $f$ with the point mass at $t = 0$, Fourier transform theory
tells us that the same value is obtained by pairing their Fourier transforms.
The Fourier transform of the point mass is the constant function $1$ but
because we know $\hat{f}$ is supported in $[- \pi, \pi]$ we can replace $1$
with $1 \cdot \chi_{[- \pi, \pi]}$, the characteristic function $[- \pi,
\pi]$. $1 \cdot \chi_{[- \pi, \pi]}$, is the Fourier transform of some
function $k_0$ in $PW$ and this discussion suggests, correctly that that
function is $k_{0,}$ the $PW$ reproducing kernel for evaluating at $t = 0$:

\begin{align*}
  k_0 (t) & = (1 \cdot \chi_{[- \pi, \pi]})^{\vee} = \frac{\sin \pi t}{\pi t}
  = \mathrm{sinc} t\\
  f (0) & = \langle f, k_0 \rangle  \text{ all } f \in PW
\end{align*}

Here \tmrsup{$\vee$} is the inverse Fourier transform, the second equality on
the first line is an elementary Fourier transform computation and the third is
the definition of the function $\mathrm{sinc}$.

This gives $k_0$, the reproducing kernel for evaluating at the origin. By
translation invariance $k_x$, the reproducing kernel for evaluating at $x$ is
$k_x (t) - \mathrm{sinc} (t - x)$ and its Fourier transform is $(k_x)^{\wedge}
(\xi) = e^{2 \pi \xi} \chi_{[- \pi, \pi]} (\xi)$. In particular the functions
$\{ (k_n)^{\wedge} \}_{n \in \mathbb{Z}}$ are an orthonormal basis of the
space of Fourier transforms of functions in $PW$. Performing the inverse
Fourier transform we see that $\{ k_n \}_{n \in \mathbb{Z}}$ is an orthnormal
basis of $PW$. Hence we have

\begin{theorem}[Shannon sampling theorem]
  If $f (t)$ is a finite energy band limited signal with spectrum contained in
  $[- \pi, \pi]$ then $f \in PW$ (by definition) and
  \begin{enumerate}
    \item the sequence of sample values $\{ f (n) \}$ is a square summable
    sequence and,
    
    \item $f$ can be reconstructed from those values using the formula
    \begin{equation}
      f (t) = \sum \langle f, k_n \rangle k_n (t) = \sum f (n) \mathrm{sinc}
      (t - n) . \label{formula}
    \end{equation}
    \item Conversely given any square summable sequence $\{ a_n \}$ there is a
    function $f$ in $PW$ with $f (n) = a_n$, for all $n$ and the value of $f$
    at all points is given by (\ref{formula}).
  \end{enumerate}
\end{theorem}

(The previous result has many names, we retreat behind the Wikipedia entry on
Stigler's law.)

This result describes the type of values obtained by regular sampling of the
function $f$ and gives a scheme for reconstructing $f$ from those sample
values---think of electronic device which samples audio signal at rate of 100
kHz and then a device which reconstructs the signal from the sample
data---think about digital music.

More generally, the space $PW$ and variations provide the mathematical
framework in which to study sampling and reconstruction of band limited
signals.

\subsection{Dirichlet space}

In this section we compare the answers to some questions for the Hardy space
with the answers to the analogous questions for the closely related Dirichlet
space. \ Some answers are very similar, some are not. Each space has a story
of its own, and we consider the Dirichlet space because much is known about it
and also, on the contrary, much is still open. We will see that sometimes the
same object of the Hardy theory has, like in a broken mirror, more than one
analog in Dirichlet theory.

The Dirichlet space $\mathcal{D}$ is the Hilbert space of holomorphic
functions on the disk. $f (z) = \sum_{n = 0}^{\infty} a_n z^n$ is in
$\mathcal{D}$ exactly if, with $\alpha = 1$, the following norm is finite:
\[ (\ast) \| f \|_{\mathcal{D}}^2 = \sum_{n = 0}^{\infty} (n + 1)^a |a_n |^2 =
   | f (0) |^2 + \frac{1}{\pi}  \int \int_{\mathbb{D}} | f' (z) |^2  (1 - | z
   |^2)^{1 - \alpha} dxdy, \]
We wrote the norm in this form to emphasize the analogy with the Hardy space
in which the formula for the norm is the case $\alpha = 0$ of the previous
formula. The parameter $\alpha$ in (*) helps highlight the close relationship
with the Hardy space. With $\alpha = 0$ the formula describes the Hardy space
norm. (With $\alpha = - 1$ that formula defines the norm of the Bergman space,
another much studied RKHS). The Dirichlet space is an RKHS and it is not hard
to verify that the reproducing kernel is
\[ k_z (w) = \frac{1}{\bar{z} w} \log \left( \frac{1}{1 - \bar{z} w} \right) .
\]
These kernel functions, as well as the kernel functions for the Hardy space,
have the property that the region where $k_z$ is relatively large is roughly
the region between $z$ and the unit circle. More specifically, if $z = re^{i
\theta}$ then the region where $k_z$ is large is, roughly, $T_z$, the
intersection of the unit disk with a disk centered at $e^{i \theta}$ of radius
$2 (1 - r)$. In particular the boundary value function has its mass
concentrated near a particular point with a specific scale of dispersion. \ We
will discuss the two parameter phase space described by position and scale
further in the next section.

The Dirichlet space has not so far found a place in signal theory. We discuss
it here because it helps illuminate the Hardy space, and, truth be told,
because the authors are very fond of it.

\subsection{The Shift operator and invariant subspaces}

The operator $M_z$ of multiplication by $z$ acts boundedly on $\mathcal{D}$.
This operator, called the Dirichlet shift, has the same action on the sequence
of Taylor coefficients of a function as the Hardy space shift does for Hardy
space functions, it shifts each entry of that sequence one place to the right.
The shift on the Hardy space is isometric and that is the starting point of an
analysis which eventual leads to the theory of inner-outer factorization of
functions and a characterization of the invariant subspaces of the shift
operator acting on $H^2$. The analysis of the invariant subspaces of the shift
operator on $\mathcal{D}$ is more complicated and less complete than for
$H^2$.

The Dirichlet shift is bounded and it is easy to see that it has lots of
invariant subspaces. In particular the structure of the invariant subspaces of
finite codimension is exactly the same as for $H^2$; they are the subspaces of
functions which vanish on a given finite point sets. Some other properties of
the shift invariant subspaces of $H^2$ which follow easily from Beurling's
theorem are also true for spaces invariant of the Dirichlet shift, but with
proofs that are less straightforward and more subtle. Two examples are the
fact that any invariant subspace contains a bounded function and the fact that
the intersection of any two nontrivial invariant subspace contains a third.

There is not yet a description of the shift invariant subspaces of
$\mathcal{D}$. In fact it is not yet known how to characterize the functions
with the property that the smallest closed invariant subspace containing them
is the whole space. For the Hardy space those functions are exactly the outer
functions. For the Dirichlet space the functions must be Hardy space outer
functions and the set on which their boundary values are zero must be a
Dirichlet space null set. (A Hardy space function is the zero function if its
boundary values are zero on a set of positive Lebesgue measure. The analogous
statement for the Dirichlet space holds for smaller sets, those of logarithmic
capacity zero.) It was conjectured by Brown and Shields in 1984 that those two
conditions characterize the Dirichlet space analogs of outer functions.

\subsubsection{Multiplication operators, Carleson measures, Hankel forms}

Multiplication by the coordinate function is a bounded operator on
$\mathcal{D}$ and it follows that multiplication by a polynomial is a bounded
operator on $\mathcal{D}$. It is then natural to ask what are the multipliers
of $\mathcal{D}$, the functions $b$ such that $M_b$, multiplication by $b$ is
a bounded map of $\mathcal{D}$ into itself. (Elements in a RKHS are functions
on a set and hence there is a natural way to multiply two of them. The
question of characterizing the multipliers makes sense on any RKHS.)

If $M_b$ is a bounded multiplier on the Dirichlet space then $b$ must be a
bounded function; in fact the argument is the same as for the Hardy space
multipliers, an argument that works for any RKHS. Also $b$ must be
holomorphic. Those conditions, $b \in H^{\infty}$, are the full story for the
Hardy space but not for the Dirichlet space. To see why not select $f \in
\mathcal{D}$ and consider the requirement that $bf \in \mathcal{D}$. By
definition we must have that $(bf)' = b' f + bf'$ is square integrable.
Because $f \in \mathcal{D}$ and $b$ is bounded the second term is. Requiring
the first term to be square integable, for every $f \in \mathcal{D}$, leads to
the definition of Carleson measure for $\mathcal{D}$.

A measure $\mu$ on $\mathbb{D}$ is a Carleson measure for $\mathcal{D}$ if
\[ [\mu]_{\mathrm{CM} (\mathcal{D})} = \sup_{f \in \mathcal{D}}  \frac{\int
   \int |f|^2 d \mu}{\|f\|_{\mathcal{D}}^2} = \left\| \mathrm{Id}
   \right\|_{\mathcal{B} (\mathcal{D}, L^2 (\mu))}^2 < \infty . \]
We define $X$ to be the space of holomorphic functions $b$ defined on the disk
such that $| b' |^2 dxdy$ is a Carleson measure for $\mathcal{D}$. Considering
$f = 1$ in the previous definition we see that $X \subset \mathcal{D}$.

Our analysis to this point shows that if $M_b$ is a bounded multiplication
operator then $b \in X \cap H^{\infty}$. The argument is easily reversed and
we have the full story.

\begin{theorem}
  $M_b$ is a bounded multiplication operator on $\mathcal{D}$ if and only if
  $b \in X \cap H^{\infty}$,
\end{theorem}

Although this does not look like our description of bounded multiplication
operators for the Hardy space, it is in fact very similar. Using the
description of the Hardy space given by by (*) with $\alpha = 0$ and then
following the ideas in that section will lead to the conclusion that $M_b$ is
a bounded multiplication operator on $H^2$ if and only if $b \in BMO \cap
H^{\infty}$, which is the analog of the previous theorem. However that last
statement can be simplified because $H^{\infty} \subset BMO$, the analogous
simplification is not possible for the Dirichlet space because $H^{\infty}
\nsubseteq X$.

Of course our understanding of the space $X$ is limited by how well we
understand Dirichlet space Carleson measures. There are several known
characterizations of those measures, some are measure theoretic "local $T 1$
conditions" others are in terms of logarithmic capacity. The appearence of
logarithmic capacity does not come as a surprise: functions in $\mathcal{D}$
are defined by a Sobolev norm, and capacity has a role in the study of Sobolev
spaces somewhat similar to the role of measure theory in studying Lebesgue
spaces. However even with those results the space $X$ and the Dirichlet space
Carleson measures are much less well understood then their more classical
cousins; $BMO$ and "classical" Carleson measures.

\subsubsection{The Pick property}

Having gone this far with our analysis of multipliers for the Dirichlet space
we can consider the analog of Pick's question: Given a finite set of points in
the disk what are the necessary and sufficient conditions on a set of target
values which insure that there is a Dirichlet space multiplier of norm at most
one which takes the target values at the points of the given set.

When we looked at the similar question in the Hardy space we started by
showing that the kernel functions were eigenfunctions of the operator
$M_b^{\ast}$, the adjoint of $M_b$, and the associated eigenvalues were the
conjugates of the values of the multiplier at the given point set. This was
enough to generate a condition involving finite matrices which was necessary
in order for there to be a multiplier of the desired sort. That argument holds
for any RKHS and the matrix produced this way is called the \tmtextit{Pick
matrix} of the problem. Pick's theorem was that in the Hardy space the
condition on the Pick matrix was also sufficient for a solution to the
interpolation problem. It is now understood that there is a class of RKHS for
which an analog of Pick's theorem holds as well as a matricial version, spaces
with the \tmtextit{complete Pick property. }In recent decades it has become
clear that those RKHS have a very rich additional structure. One of the
reasons for recent interest in the Dirichlet space is that it is one of
simplest spaces other than the Hardy space with this fundamental property.

\subsubsection{Hankel forms}

On the Hardy space we considered the following bilinear Hankel form. Select a
holomorphic \tmtextit{symbol function} $b$ and define the Hankel form on the
Hardy space with symbol $b$ to be the bilinear form on $H^2$ given by, for $f,
g \in H^2$
\[ H_b  (f, g) = \langle fg, b \rangle_{H^2} . \]
We can define a Hankel form on the Dirichlet space for $f, g \in \mathcal{D}$
using the same formula but, of course, with the $\mathcal{D}$ inner product.

When we looked at Hankel forms on the Hardy space it was straightforward to
see that if $| b' |^2 dxdy$ was a Hardy space Carleson measure then $H_b$ was
bounded on the Hardy space. It then follows that having $b$ in $BMO$ will be a
sufficient condition for boundedness. The same analysis shows that having $b$
in $X$ is sufficient for $H_b$ to be bounded on the Dirichlet space. In fact,
as with the Hardy space, that is the full story.

\begin{theorem}
  The Dirichlet space Hankel form $H_b$ is bounded if and only if $b \in X$
\end{theorem}

(The definition of Hankel operators and forms for the Dirichlet space is a
place where there is more than one natural extension of the Hardy space ideas.
Emphasizing different analogies between the Dirichlet space and Hardy space
can lead to the conjugate linear map from $\mathcal{D}$ to itself given by
\[ f \rightarrow \int P (b'  \bar{f}) =\mathcal{H}_b f \]
as the natural generalization of Hankel operators to the Dirichlet space.
(Here $P$ is the orthogonal projection associated with the Bergman space.) The
condition $b \in X$ is also necessary and sufficient for is sufficient for the
boundedness of $\mathcal{H}_b$ and the proof of the easy half of the result is
the same as for $H_b$. However the full proof is different.)

The proof of the Hardy space version of the previous theorem exploited the
fact that every function in $H^1$ is the product of two functions in $H^2$ and
the duality between $H^1$ and $BMO$. Starting with the previous theorem one
can try to reverse those arguments to find our what the space $X$ is the dual
of. That leads to the notion of weakly factored spaces. We define the weakly
factored space $\mathcal{D} \odot \mathcal{D}$ to be the space of those $f$
holomorphic on $\mathbb{D}$ for which
\[ \|f\|_{\mathcal{D} \odot \mathcal{D}} = \inf \left\{ \sum_j \|g_j
   \|_{\mathcal{D}} \|h_j \|_{\mathcal{D}} : \sum_j g_j h_j = f \right\} <
   \infty . \]
A consequence of the previous theorem is the duality relation

\begin{corollary}
  $(\mathcal{D} \odot \mathcal{D})^{\ast} = X$
\end{corollary}

Using the factorization of $H^1$ functions described in Lemma \ref{pluto} it
is straightforward to see that $H^1 = H^2 \odot H^2$. Hence the previous
corollary is the Dirichlet space analog of Fefferman's classical $(H^1)^{\ast}
= BMO$.

Using interpolation of Banach spaces, real or complex, it is possible to start
from the spaces $H^1$ and $BMO$ and recover the full range of Hardy spaces
$H^p$, $1 < p < \infty$ with the starting Hilbert space $H^2$ in the middle of
the scale. Similarly one can construct the scale of spaces connecting
$\mathcal{D} \odot \mathcal{D}$ and $X$ which has the Hilbert space
$\mathcal{D}$ in the middle. Very little is known about those spaces.

\subsection{Dyadic Dirichlet Space}

Let $\mathcal{T}$ be the vertex set the dyadic tree, which we choose to also
call $\mathcal{T}$. Thus $\mathcal{T}$ is a connected, simply connected,
rooted graph with two edges at the root vertex $o$ and three edges at all the
other vertices. We put a partial order, $\preceq$, on the vertices by saying
$\alpha \preceq \beta$ exactly if $\alpha$ is a vertex on the geodesic path
connecting $o$ and $\beta$. For any $\beta \in \mathcal{T} \setminus \{ o \}$
we let $\beta^-$ be its predecessor, the maximal $\alpha$ such that $\alpha
\preceq \beta$ and $\alpha \neq \beta$.

We use two functions, $I$ and $\Delta$ acting on functions defined on
$\mathcal{T}$:

\begin{align*}
  If (\beta) & = \sum_{o \preceq \alpha \preceq \beta} f (\alpha),\\
  \Delta f (\beta) & = \left\{ \begin{array}{cc}
    0 & \text{if } \beta = o\\
    f (\beta) - f (\beta^-) & \text{otherwise}
  \end{array} \right. .
\end{align*}

These operators are models for integration and ifferentiation. If $f$ is a
function on $\mathcal{T}$ with $f (o) = 0$ then $I \Delta f = \Delta If = f$.
We define the dyadic Dirichlet space, $\mathcal{D}_{\mathtt{dyad}}$, to be the
Hilbert space of functions $f$ defined on $\mathcal{T}$ for which $\Delta f
\in \ell^2 (\mathcal{T)}$. The space is normed by
\[ \| f \|_{\mathcal{D}_{\mathtt{dyad}}}^2 = | f (o) |^2 + \| \Delta f
   \|_{\ell^2 (\mathcal{T)}}^2 . \]
This space is a RKHS, the reproducing kernel for evaluation at $\alpha \in
\mathcal{T}$ is $k_{\alpha} = I (\chi_{[o, \alpha]})$.

\subsubsection{$\mathcal{D}_{\mathtt{dyad}}$ is a model for $\mathcal{D}$}

One of the reasons for considering the space $\mathcal{D}_{\mathtt{dyad}}$ is
that it is a simple model for $\mathcal{D}$. The analogy is best understood by
regarding $\mathcal{T}$ as a point set in the unit disk. Informally, the root
is placed at the origin, the $2^n$ vertices connected to the origin by
geodesics of length $n$ are spaced evenly on the circle of radius $1 - 2^{-
n}$. The edge between an $\alpha$ on that circle to its predecessor $\alpha^-$
is represented by an almost radial line segment connecting the two.

In this picture the values of an $f \in \mathcal{D}_{\mathtt{dyad}}$ at points
of the abstract tree are a model for the values of some unspecified function
$\tilde{f}$ $\in \mathcal{D}$. If fact starting with any $g \in \mathcal{D}$
and restricting to the points of the realization of $\mathcal{T}$ inside the
disk will given an element of $\mathcal{D}_{\mathtt{dyad}}$. Continuing the
analogy, if $f \in \mathcal{D}_{\mathtt{dyad}}$ then $\Delta f$ is a model for
$\tilde{f}'$ and the fact that $\Delta f$ is required to be square summable
models the fact that $\tilde{f}'$ must be square integrable. (Our view from
great height is ignoring scaling: $\Delta f$ is actually a model of the
invariant derivative $\delta f (z) = (1 - | z |^2) f' (z)$.)

\subsubsection{The results are similar}

The analogies just described are relatively superficial. More interesting is
that the analogies extend to subtle aspects of the Dirichlet space theory.
There are natural extensions of the definitions of multipliers, of Carleson
measures, of Hankel forms, etc. from the Dirichlet space to the dyadic
Dirichlet space. For all of the results we have discussed (and many others)
the results for the two spaces are "the same", that is they continue the
pattern suggested by the analogy. Generally the proofs in the dyadic case are
easier and sometimes those proofs provide road maps for the more difficult
proofs for the classical space.

Carleson measures are a particularly interesting case. The measure theoretic
characterization of Carleson measures for $\mathcal{D}$ is most simply
obtained by first solving the analogous problem in
$\mathcal{D}_{\mathtt{dyad}}$ and then using the fact mentioned before, that
the restriction of functions in $\mathcal{D}$ produces functions in
$\mathcal{D}_{\mathtt{dyad}}$, to lift the result to $\mathcal{D}$.

\subsubsection{Phase space analysis}

A number of interrelated ideas form the general category of phase space
analysis. The RKHS we have discussed are in this category and the dyadic
Dirichlet space is a particularly simply instance. We will say a few words
about the general theme but, even by the standards of what has gone before, we
will be very informal. Our main point is that some of the ideas we have seen
here are instances of general themes.

Suppose we wanted to analyze a function $f$ in the Dirichlet space. We know
there are reproducing kernels $\{ k_z \}_{z \in \mathbb{D}}$ and hence $f (z)
= \langle f, k_z \rangle$. We mentioned that reproducing kernels were a
substitute for an orthonormal basis. If they were an orthogonal basis we would
have a representation
\begin{equation}
  f = \sum \left\langle \dot{f}, \frac{k_z}{\| k_z \|} \right\rangle
  \frac{k_z}{\| k_z \|} \label{A}
\end{equation}
but that is not true. A possible path forward is to replace the sum by an
integral and hope for a representation
\begin{equation}
  f = \int \langle \dot{f}, k_z \rangle k_z d \mu (z) . \label{B}
\end{equation}
Here we have absorbed the normalizing factors into the measure but we are
intentionally vague about the details. This does not hold but a formula of
this type is true for the Bergman space ("Bergman reproducing formula") and in
a number of spaces of interest in quantum theory ("coherent state
representations"). Another way to try to go forward is to try to use a subset
of the $\{ k_x \}$ and obtain a summation formula of the type (\ref{A}), for
instance using only those $z$ which correspond to the vertices of
$\mathcal{T}$. That set is still not an orthogonal basis but it is close
enough so that (\ref{A}), while not latterly true, is a good enough
approximation, both analytically and conceptually, to be a useful starting
point. That fact is the heart of the relation between
$\mathcal{D}_{\mathtt{dyad}}$ and $\mathcal{D}$. It is also the starting point
for obtaining representations of functions in various function spaces as
linear combinations of reproducing kernels associated with points in a set
such as $\mathcal{T}$.

When we discussed the Hardy space there were different viewpoints; Hardy space
functions can be viewed as holomorphic functions in the disk or as boundary
value functions on the circle, and it is possible to pass back and forth
between those viewpoints with no loss of information. The same is true for
many other spaces of functions on the disk. Consider now how that interacts
with formulas such as (\ref{A}) and (\ref{B}) and their various refinements.
We could start with a boundary function $f_{\mathtt{boundary}}$ pass to the
associated function inside the disk, $f_{\mathtt{inside}}$, use the analytical
tools to represent $f_{\mathtt{inside}}$ as a sum or integral of simple
pieces, and then pass back to the boundary function. This would realize
$f_{\mathtt{boundary}}$ as a sum (or integral) of boundary values of a set of
well understood functions. If the coefficient corresponding to $z$ in the
representations is built by taking the inner product of $f_{\mathtt{inside}}$
with some $h_z$, function concentrated on the set we introduced earlier,
$T_z$, then it will be mainly responsive to the values of
$f_{\mathtt{inside}}$ inside $T_z$ and hence presumably to the values of
$f_{\mathtt{boundary}}$ near the part of the unit circle cut off by $T_z$.
Furthermore the boundary values of the function in the representation, perhaps
again $h_z$, will also be concentrated on that same interval. In sum, the
representation of a function on the boundary uses analysis and reconstruction
tools paramertrized by two real parameters. The parameters can be understood
as position and scale, the center of the boundary interval and its length, and
those parameters form points in "phase space". For the Hardy space the points
$z = re^{i \theta}$ parameterize the disk $\mathbb{D}$ which is the phase
space; $re^{i \theta}$ is the complex parameter describing the interval on the
circle with center $e^{i \theta}$ and radius $1 - r$.

Without examples the previous paragraph is idle talk. However there are
examples. Many RKHS of holomorphic functions in one and several complex
variables fit this pattern, or they do after minor modifications. The Bergman
spaces are fundamental examples. Also there is an important class of examples
not related to holomorphic functions. It is possible to start with a general
function on the circle, or on the line, or on $n$-space and form an associated
phase space, a space of one higher dimension whose new coordinate is scale.
There are systematic ways to extend a function $f$ on the space to a function
$f_{\mathtt{inside}}$ defined on the phase space and to introduce functions
$\{ k_{\zeta} \}$ for $\zeta$ in the phase space. and proceed exactly as
described. With the appropriate details filled in the result is an exact
formula in the style of (\ref{B}). The functions $k_{\zeta}$ are each
associated with a point in phase space and their boundary values, their traces
on the starting space, are concentrated in the associated ball, the ball whose
center and radius are the coordinates in phase space. In fact all this can be
done with the $k_{\zeta}$ all translates and dilates of a single function, a
"mother wavelet". The resulting formula is the "Calderon reproducing formula"
or the "continuous wavelet transform". There is a striking refinement of these
ideas. It is possible to arrange the details so that the set of normalized
$k_{\zeta}$ with $\zeta$ in a discrete subset of phase space, shaped like
$\mathcal{T}$, is an orthonormal basis of the Lebesgue space of the starting
manifold. In that case there is a discrete representation, a formula of the
form (\ref{A}) for representing any function. In that formula the coefficients
and the summands, the analysis and the reconstruction, respect the description
of the function in terms of the phase space parameters of location and scale.
The resulting formula is the "wavelet representation" of the function which is
fundamental in large areas of signal analysis.

\section{Further reading}

\begin{itemize}
  \item A lovely and quick introduction to some of the topics we have discussed
  is the self-contained, expository article
  
  \
  
  John McCarthy Pick's theorem - what's the big deal?~{\tmem{American
  Mathematical Monthly~}}Vol. 110 No. 1 [2003] 36-45,
  
  \
  
  {\noindent}where in a few pages the route from the Hardy space to control
  theory to Pick's theory is covered.
  
  \item {\noindent}The pure mathematician who wants to painlessely understand
  what signal theory and the related control theory are about, can watch the
  old, but clear and enjoyable, 1987 MIT lectures of Alan Oppenheim,
  
  \
  
  https://ocw.mit.edu/resources/res-6-007-signals-and-systems-spring-2011/video-lectures/
  
  \
  
  {\noindent}where some surprisingly effective pratical applications are
  shown.
  
  \item A very nice introduction to {\noindent}$H^{\infty}$ control theory are
  the 2008 lecture notes for ``the mythical 'mathematically mature engineering
  student'' at University of Toronto,
  
\

  https://www.control.utoronto.ca/\~{}broucke/ece356s/ece356Book2008.pdf
  
  \
  
  {\noindent}by Bruce Francis, one of the protagonists of contemporary
  holomorphic control theory.
  
  \item A largely overlapping body of knowledge, but from the viewpoint of the
  pure mathematician, is in the monograph
  
  \
  
  Jonathan R. Partington - Linear operators and linear systems: An analytical
  approach to control theory (2004, CUP)
  
  \
  
  {\noindent}which also works as a comprehensive introduction to Hardy space
  theory.
  
  \item {\noindent}An excellent survey (with proofs) on Hankel operators and
  Nehari theory is
  
  \
  
  Vladimir Peller, An Excursion into the Theory of Hankel Operators,
  Holomorphic Spaces MSRI Publications Volume 33, 1998,
  
  \
  
  which can be found here: http://mathscinet.ru/files/PellerV.pdf
  
  \item An excellent, self-contained, and easy to read monograph on
  reproducing kernel Hilbert spaces and Pick theory, also providing an
  introduction to Hardy space theory, is
  
  \
  
  Jim Agler, John McCarthy, Pick Interpolation and Hilbert Function Spaces,
  American Mathematical Society, 2002.
  
  \item The discourse on Nehari, Hankel, Toeplitz, Hilbert transform, and BMO,
  is the subject ofthe short and dense
  
  \
  
  Donald Sarason, Function Theory on the Unit Disc, Virginia Polytechnic
  Institute and State University, 1978
  
  \item To move deeper in hard-analysis Hardy space theory, our standard
  reference is still
  
  \
  
  John Garnett, Bounded analytic functions, Springer, Revised 1st ed. 2007
  
  \item A standard text of Functional Analysis which is fully adequate for the
  subject is
  
  \
  
  Peter Lax, Functional Analysis, Wiley 2002. \
  
  \item A chapter on the Paley-Wiener space, with a thourogh discussion of
  sampling results (which are crucial in applications to engineering) is
  
  \
  
  Kristian Seip, Interpolation and Sampling in Spaces of Analytic Functions,
  American Mathematical Soc., 2004.
  
  \item There are two recent monographs on the Dirichlet space:
  
  \
  
  Omar El-Fallah, Karim Kellay, Javad Mashreghi, Thomas Ransford, A primer on
  the Dirichlet Space, Cambridge Tracts in Mathematics, 2014,
  
  \
  
  and
  
  \
  
  Nicola Arcozzi, Richard Rochberg, Eric T. Sawyer, Brett D. Wick, The
  Dirichlet Space and Related Function Spaces, American Mathematical Society,
  2019.
  
  \
  
  The former develops the theory from a classical point of view, the latter
  from the viewpoint of Reproducing Kernel Hilbert Spaces.
  
  \item  An excellent way to become acquainted to time-frequency analysis is
  
  \
  
  Ingrid Daubechies, Ten lectures on Wavelets, SIAM, 1994,
  
  \
  
  by one of the pioneers of wavelet theory.
  
  \item Specific operators on specific Hilbert function spaces can ``model''
  general classes of operators acting on Hilbert spaces. This line of
  investigation has one of its milestones in:
  
  \
  
  B. Sz. Nagy and C. Foias, Harmonic Analysis of Operators on Hilbert Space.
  VIII + 387 S. Budapest/Amsterdam/London 1970. North Holland Publishing
  Company
  
  \item Finally, we suggest this classical, short monograph, where the ideas
  surrounding Beurling's theorem on invariant subspaces are the starting point
  to derive in a simple way some deep results in Hardy space theory:
  
  \
  
  Helson, Henry Lectures on invariant subspaces. Academic Press, New
  York-London 1964 xi+130 pp
\end{itemize}

\end{document}